\documentclass[12pt,a4paper]{article}
\usepackage[latin1]{inputenc} 
\usepackage[leqno]{amsmath}
\usepackage{amsthm} 
\usepackage{amssymb} 
\usepackage{amstext} 
\usepackage{amsxtra}
\usepackage{esint} 
\usepackage{mathrsfs} 
\usepackage{todonotes} 
\usepackage{cite} 
\usepackage{hyperref}
\usepackage{graphicx}
\usepackage{caption}
\usepackage{fancyhdr}
\title{Global Dynamics, Blow-Up, \\
and Bianchi Cosmology}
\author{Nitsan Ben-Gal, Bernhard Brehm, Johannes Buchner, \\
Juliette Hell, Anna Karnauhova, Stefan Liebscher,\\ 
Alan Rendall, Brian Smith, Hannes Stuke,\\
Martin V\"ath, Bernold Fiedler\\
\vspace{2cm}}

\date{version of \today}

\begin{document}
\maketitle
\thispagestyle{empty}

\vfill

Institut f\"ur Mathematik\\
Freie Universit\"at Berlin\\
Arnimallee 7\\ 
14195 Berlin, Germany


\newpage
\pagestyle{plain}
\pagenumbering{roman}
\setcounter{page}{1}

\begin{abstract}
Many central problems in geometry, topology, and mathematical physics lead to questions concerning the long-time dynamics of solutions to ordinary and partial differential equations. Examples range from the Einstein field equations of general relativity to quasilinear reaction-advection-diffusion equations of parabolic type. Specific questions concern the convergence to equilibria, the existence of periodic, homoclinic, and heteroclinic solutions, and the existence and geometric structure of global attractors. On the other hand, many solutions develop singularities in finite time. The singularities have to be analyzed in detail before attempting to extend solutions beyond their singularities, or to understand their geometry in conjunction with globally bounded solutions. In this context we have also aimed at global qualitative descriptions of blow-up and grow-up phenomena. \\

We survey some of our results supported in the framework of the Sonderforschungsbereich 647 ``Space -- Time -- Matter'' of the Deutsche Forschungsgemeinschaft.
\end{abstract}


\newpage
\pagenumbering{arabic}
\setcounter{page}{1}

\section{Introduction}\label{sec:intro}

Many central problems in geometry, topology, and mathematical physics lead to questions concerning the long-time dynamics of solutions to ordinary and parabolic partial differential equations. Three basic scenarios are possible. \\

First, under suitable dissipativity conditions, longtime existence is guaranteed and a low-dimensional global attractor of uniformly bounded eternal solutions exists. The structure of such global attractors details longtime behavior. Some of those structures are encoded in meanders -- which open doors to a broad variety of algebraic structures far beyond the original PDE scope of the attractor question. \\ 

Second, longtime existence without dissipativity features trajectories which escape to infinity in infinite time (grow-up). In some cases it is possible to construct compactifications of an unbounded attractor containing, for example, equilibria at infinity. As an example we show how equivariant bifurcations give rise to black hole initial data that is not spherically symmetric. \\

Third, singularities may develop after finite time. In other words, some solutions escape to infinity within finite time (blow-up). Some of the ramifications are explored in complex time, here. Cosmologically, the most meaningful singularity is the big-bang. We address the Bianchi models: a system of ordinary differential equations which represent the homogenous, but anisotropic, Einstein equations. The Belinsky--Khalatnikov--Lifshitz (BKL) conjecture favors them as models for a spacelike singularity. We obtain stability results for a primordial universe tumbling between different Kasner states, each of which represents a self-similarly expanding universe. Furthermore we address the BKL question of particle horizons. To conclude, we embed general relativity in the broader context of Ho\v{r}ava--Lifschitz gravity, to study the influence on the tumbling behavior between the Kasner states and the resulting chaoticity. 

\section{Global attractors of parabolic equations}\label{sec:globalattractors}

\numberwithin{equation}{section}
\numberwithin{figure}{section}

In this section we survey recent progress on the global dynamics, mostly of the scalar parabolic semilinear equation
\begin{equation}\label{C8:eq:2.1.1}
u_t = u_{xx} +f\,.
\end{equation} 

The emerging global theory gives rise to surprising topological structures, as finite regular cell complexes, in Section~\ref{subsec:sturmglobalattractors}, and to combinatorial algebraic structures like Temperley--Lieb algebras, in Section~\ref{subsec:meanders}. An application to axisymmetric, but not spherically symmetric, initial data for black holes which satisfy the Einstein constraint is given in Section~\ref{subsec:blackholedata}. The question of stability of stationary solutions is addressed in Section~\ref{subsec:statsol}. \\

For simplicity we consider the PDE~\eqref{C8:eq:2.1.1} under Neumann $\mathcal{N}$ or periodic $\mathcal{P}$ boundary conditions, on the unit interval $0<x<1$. Other separated linear boundary conditions can be treated analogously to the Neumann case. The quasilinear uniformly parabolic case may easily be believed to yield identical results, but the long necessary groundwork is not complete at this stage.

\subsection{Sturm global attractors}\label{subsec:sturmglobalattractors}

We consider the PDE~\eqref{C8:eq:2.1.1} for several classes of nonlinearities $f$. In particular, we survey our results~\cite{FiedlerRochaWolfrum2012-SturmAttractors, Fiedler15, FiedlerRochaWolfrum12, FiedlerGrottaRocha14, RochaFiedler14, FiedlerRocha14}. We denote these classes by $\mathcal{N}(x,u,u_x)$, $\mathcal{N}(u,u_x)$, $\dots$ for $f=f(x,u,u_x)$, $f=f(u,u_x)$, $\dots$ under Neumann boundary conditions, and by $\mathcal{P}(x,u,u_x)$, $\mathcal{P}(u,u_x)$, $\dots$ in the $x$-periodic case. Dirichlet or other separated boundary conditions of mixed type can be treated analogously to the Neumann case. We always assume $f\in C^1$ to be dissipative, i.e., the set of bounded \emph{eternal solutions} $u(t,x),t\in\mathbb{R}$ is itself bounded in a suitable Sobolev type interpolation space $u(t,\cdot)\in X \hookrightarrow C^1$. Explicit assumptions on $f$ which guarantee dissipativity involve, for example, sign conditions $f(x,u,u_x) \cdot u<0$, for large $|u|$, and subquadratic growth in $u_x$. The bounded set $\mathcal{A}_f \subseteq X$ of spatial profiles $u(t,\cdot)$ on bounded eternal solutions is called the \emph{global attractor} of~\eqref{C8:eq:2.1.1}. The concept originated with Ladyzhenskaya's work on the two-dimensional Navier--Stokes equation. \\

For general surveys on global attractors $\mathcal{A}$ see for example~\cite{hale,vishikchepyzhov}. The general theory establishes properties like existence, global attractivity (hence the name), compactness, finite-dimensionality in the fractal box counting sense, and sometimes connectivity. From such a general, and hence necessarily superficial, point of view the global attractors $\mathcal{A}_f$ of our model class~\eqref{C8:eq:2.1.1} are just an easy exercise. \\

\begin{figure}
\centering \includegraphics[angle=0,width=\textwidth]
{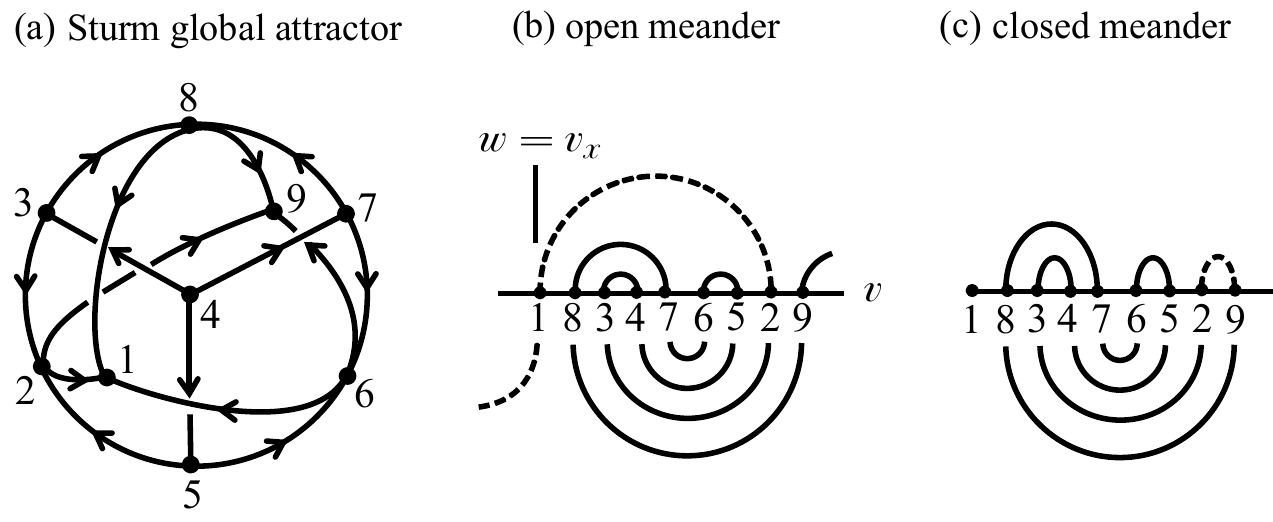}
\caption{\emph{(a) Sturm global attractor $\mathcal{A}_f$ and connection graph $\mathcal{C}_f$ of a PDE~\eqref{C8:eq:2.1.1} with 9 hyperbolic equilibria $v$, labeled by $1,\dots, 9=N$.
One-dimensional heteroclinic orbits between equilibria of adjacent Morse index are indicated by arrows.
(b) Open meander, at $x=1$, which leads to the 9 equilibrium solutions $v$ by shooting from the Neumann boundary condition at $x=0$.
(c) Closed meander which arises by replacing the dashed upper arch $(1,2)$, in (b), with the new dashed upper arch $(9,2)$; see~\cite{FiedlerCastaneda2012-Meander}.}}
\label{C8:fig:2.1.1}
\end{figure}

However, our objective is a much deeper global understanding of the detailed global structure of the global attractors $\mathcal{A}_f$. Let us address the general Neumann case $f\in\mathcal{N}(x,u,u_x)$ first. By results of Zelenyak and Matano, the PDE~\eqref{C8:eq:2.1.1} then possesses a gradient-like Morse structure due to a decreasing Lyapunov, or energy, or Morse, functional $V$. Throughout we assume all Neumann equilibrium solutions
\begin{equation}\label{C8:eq:2.1.2}
0=v_{xx}+f
\end{equation}
to be hyperbolic, i.e., to be quadratically nondegenerate critical points of $V$. Then
\begin{equation}\label{C8:eq:2.1.3}
\mathcal{A}_f =\mathcal{E}_f \cup\mathcal{H}_f
\end{equation}
consists of the finitely many equilibria $v\in\mathcal{E}_f$ and \emph{heteroclinic orbits} $u(t,\cdot)\in\mathcal{H}_f$ \emph{connecting} them:
\begin{equation}\label{C8:eq:2.1.4}
u(t,\cdot)\rightarrow v_{\pm}\in\mathcal{E}_f \quad \text{for} \quad t \rightarrow \pm \infty\,,
\end{equation}
with equilibria $v_+\neq v_-$. We call $\mathcal{A}_f$ a \emph{Sturm global attractor}, to emphasize the particular origin from the one-dimensional parabolic PDE class~\eqref{C8:eq:2.1.1}. Information on the vertex pairs $(v_-,v_+)$ which do possess a heteroclinic edge~\eqref{C8:eq:2.1.4}, and on the pairs which don't, can be encoded in the (directed, acyclic, loop-free) \emph{connection graph} $\mathcal{C}_f$; see~\cite{FiedlerCastaneda2012-Meander} and Fig.~\ref{C8:fig:2.1.1}(a). \\

As was discovered by Fusco and Rocha~\cite{FuscoRocha}, a permutation $\sigma$ given by the boundary values of the equilibria $\mathcal{E}_f= \{ v_1,\dots, v_N\}$ is the key to a global analysis of the Sturm global attractor $\mathcal{A}_f$ and its heteroclinic connection graph $\mathcal{C}_f$. The \emph{Sturm permutation} $\sigma=\sigma_f$ encodes the orderings of the equilibria,
\begin{equation}\label{C8:eq:2.1.5}
\begin{aligned}
v_1 			\;&<\; v_2 			\!\!\!\!&<\; \dots 	\;&< \; v_N			&\quad \text{at}\quad x=0\,,\\
v_{\sigma (1)} 	\;&<\; v_{\sigma (2)} 	\!\!\!\!&< \;\dots 	\;&< \;v_{\sigma (N)}	&\quad \text{at} \quad x=1\,,
\end{aligned}
\end{equation}
at the respective Neumann boundaries. More precisely, the Sturm permutations $\sigma_f$ determine all \emph{unstable dimensions}, alias \emph{Morse indices}, $i=i(v_k)$ of the equilibria $v_k \in \mathcal{E}$; see~\cite{FuscoRocha}. This is a linear ODE question of Sturm--Liouville type; hence the name Sturm permutation for $\sigma_f$. Moreover $\sigma_f$ determines the (PDE) heteroclinic connection graph $\mathcal{C}_f$, explicitly and constructively, in a combinatorial and algorithmic style~\cite{FiedlerRocha1996-SturmPDE}. Global Sturm attractors $\mathcal{A}_f$ and $\mathcal{A}_g$ with the same Sturm permutation $\sigma_f = \sigma_g$ are in fact $C^0$ orbit equivalent, i.e.\ homeomorphic in $X$ under a homeomorphism which maps eternal PDE orbits $\{ u(t,\cdot),t\in\mathbb{R}\}\subseteq\mathcal{A}_f$ to eternal orbits in $\mathcal{A}_g$; see~\cite{FiedlerRocha00}. \\

Not all permutations $\sigma$ are Sturm permutations $\sigma_f$ arising from the Neumann class $f \in \mathcal{N}(x,u,u_x)$ of~\eqref{C8:eq:2.1.1}. In fact $\sigma=\sigma_f$ for some $f$, if and only if the permutation $\sigma\in S_N$ is a meander permutation which, in addition, is dissipative and Morse~\cite{FiedlerRocha1999-SturmPermutations}. Following Arnold, here, we call $\sigma$ a \emph{meander permutation} if $\sigma$ describes the transverse crossings of a planar $C^1$ Jordan curve, or ``meandering river'', oriented from south-west to north-east, with the horizontal axis, or ``road''. Labeling the crossings along the ``river'', sequentially, the meander permutation $\sigma$ is defined by the labels $\sigma(1),\sigma(2),\dots,\sigma(N)$ read along the road, left to right. We call a meander $\sigma$ \emph{dissipative} if $\sigma(1)=1$ and $\sigma(N)=N$.
\begin{figure}
\centering\includegraphics[angle=0,width=0.80\textwidth]
{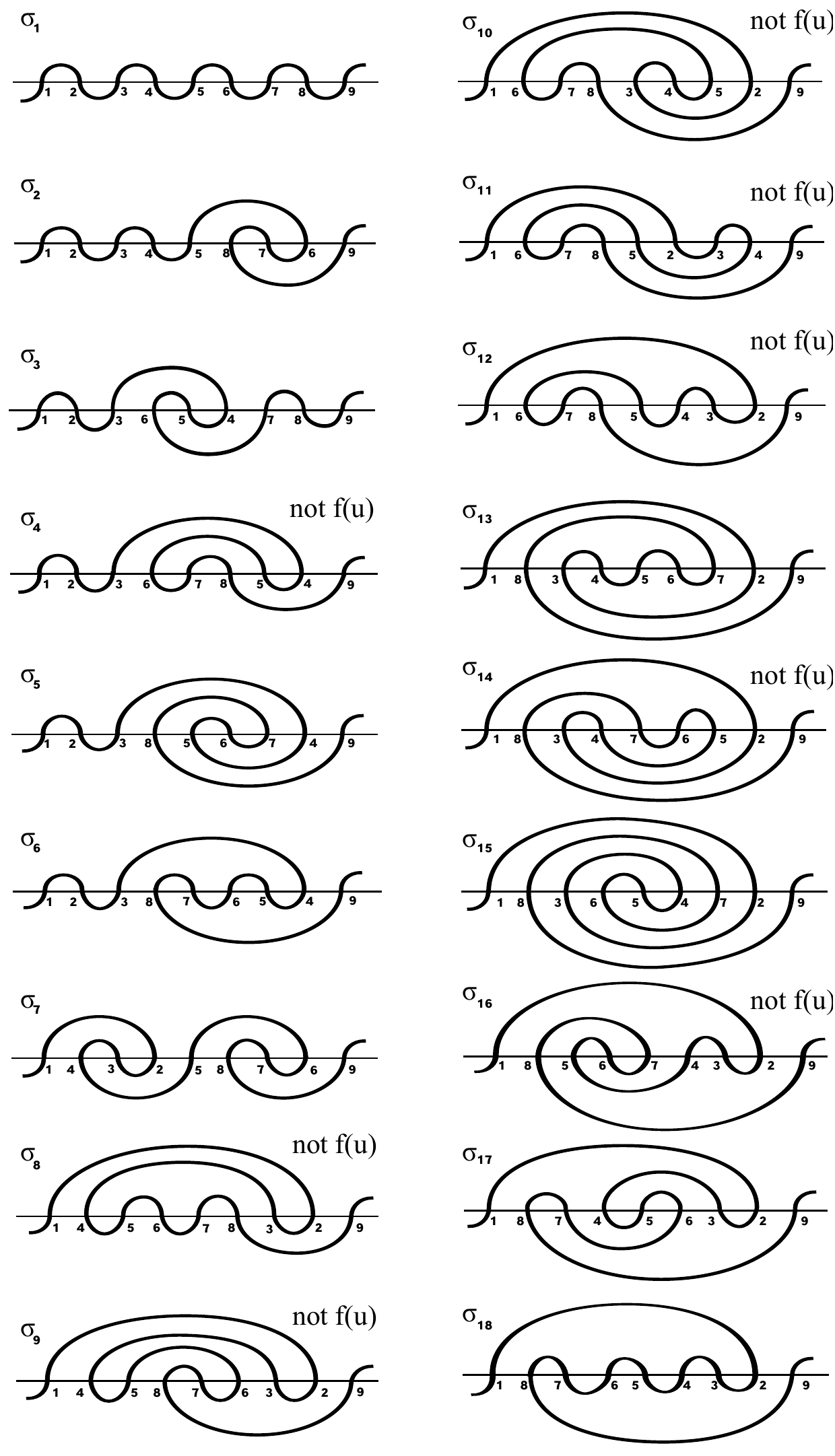}
\caption{\emph{All Sturm permutations in $S(9)$, up to trivial symmetries, and their meanders. We indicate all cases which cannot be realized by nonlinearities $f =f(u)$; see~\cite{FiedlerRochaWolfrum2012-SturmAttractors}.}}
\label{C8:fig:2.1.2}
\end{figure}
We call $\sigma$ \emph{Morse} if the integers
\begin{equation}\label{C8:eq:2.1.6}
i_k:= \sum_{j=1}^{k-1}(-1)^{j+1}\operatorname{sign}\big(\sigma^{-1}(j+1)-\sigma^{-1} (j)\big)
\end{equation}
are nonnegative, for all $k=1,\dots,N$. Geometrically, this amounts to a net right winding of the river, as counted from the first crossing; see Fig.~\ref{C8:fig:2.1.1}(b) and Fig.~\ref{C8:fig:2.1.2} for many examples. The PDE connection graphs $\mathcal{C}_f$ then follow, in each case, without further PDE analysis. \\

It is not difficult to see why Sturm permutations $\sigma=\sigma_f$ satisfy these three properties. Dissipativeness of $\sigma$ follows because $\mathcal{E}_f$ must contain extreme equilibria $v_1(x)<v_2(x),\dots,v_{N-1}(x)<v_N(x)$, for all $0\leq x\leq1$, by dissipativeness of~\eqref{C8:eq:2.1.1} and parabolic comparison. The Morse property $i_k\geq0$ follows because $i_k=i(v_k)\geq0$ are the nonnegative Morse indices; see~\cite{FuscoRocha}. The meander property follows by ODE shooting for the equilibrium boundary value problem~\eqref{C8:eq:2.1.2} in the $(v,v_x)$ plane. Indeed, the ``river'' is the image of the horizontal Neumann $v$-axis $\{ v_x =0\}$, at $x=0$, propagated to $x=1$ under the ODE-flow~\eqref{C8:eq:2.1.2}. The ``road'' is the horizontal Neumann $v$-axis $\{v_x=0\}$, at $x=1$. See Section~\ref{subsec:meanders} for a deeper discussion of meanders. The converse claim, how all dissipative Morse meanders $\sigma$ are actually Sturm permutations $\sigma=\sigma_f$ for some $f$, is beyond the limitations of this introductory survey. See \cite{FiedlerRocha1999-SturmPermutations}. The deceptively elementary description of Sturm permutations $\sigma_f$ should not conceal the formidable complications arising from the meander property. Already Gauss got intrigued by the related question of self-intersecting orders in planar curves; see~\cite{Gauss}. Rescaling $x=1$ to become large, the evident complexities of chaotic stroboscope maps for the periodically forced pendulum or the forced van~der~Pol oscillator are also lurking here, as special cases. \\

A key technical ingredient to the results indicated above are nodal properties, i.e.\ a refinement of the parabolic comparison principle in one space dimension, or in radially symmetric settings. Indeed, let $u_1(t,x)$ and $u_2(t,x)$ be any two non-identical solutions of~\eqref{C8:eq:2.1.1}, in $X \hookrightarrow C^1$, and let the \emph{zero number} $z(\varphi)\leq\infty$ count the strict sign changes of any continuous $x$-profile $\varphi: [0,1]\rightarrow\mathbb{R}$. Then the zero number
\begin{equation}\label{C8:eq:2.1.7}
t \mapsto z\big(u_1(t, \cdot) -u_2 (t,\cdot )\big)
\end{equation}
is finite, for any $t >0$, nonincreasing with $t$, and drops strictly at any multiple zeros of $x\mapsto u_1(t,x)-u_2(t,x)$; see~\cite{Angenent}. Early versions of this, in the linear case, in fact date back to Sturm and convinced us to name the nonlinear ramifications after him; see~\cite{Sturm}. \\

It is the zero number dropping property~\eqref{C8:eq:2.1.7} which generates surprising rigidity in the global Sturm attractors $\mathcal{A}_f$ and enables a detailed global analysis. For example, stable and unstable manifolds of hyperbolic equilibria $v_+$ and $v_-$ intersect transversely, automatically, along any heteroclinic PDE orbit $u(t, \cdot)\rightarrow v_\pm$, for $t \rightarrow\pm\infty$, in $\mathcal{H}_f$; see~\cite{Angenent, Henry85}. In particular, this implies local structural stability of Sturm attractors and, moreover, a strict dropping of Morse indices,
\begin{equation}\label{C8:eq:2.1.8}
i(v_-) > z\big(u(t,\cdot)-v_-\big) \geq z(v_+ -v_-) \geq z\big(v_+ -u(t,\cdot)\big) \geq i(v_+)\,.
\end{equation}
By a (standard) transitivity property of transverse heteroclinicity and a (Sturm specific) cascading principle it turns out, in fact, that the full heteroclinic connection graph $\mathcal{C}_f$ is generated by the edges $(v_-,v_+)$ with adjacent Morse index. Saddle-saddle heteroclinic orbits $i(v_-) = i(v_+)$ are also excluded by Morse dropping~\eqref{C8:eq:2.1.8}. This precludes global topological ambiguities which, for example, homotopy-invariant Conley--Morse theory is not allowed to discard. Indeed the positive integer polynomial terms $(1+t)\,q(t)$ in the Morse (in)equalities, alias connecting homomorphisms or Conley--Franzosa connection matrices, always need to accommodate global reshufflings of the connection graph due to nongeneric saddle-saddle connections. This is where nodal arguments involving the zero number $z(u_1-u_2)$ step in. \\

We summarize more recent results on Sturm global attractors $\mathcal{A}_f$ and their characterizing Sturm permutations $\sigma_f$ for the parabolic PDE~\eqref{C8:eq:2.1.1} under Neumann $(\mathcal{N})$ and periodic $(\mathcal{P})$ boundary conditions, next. We conclude with recent progress in our global geometric understanding of $\mathcal{A}_f$ as regular topological cell complexes. \\

Little information is provided, in our approach, on the global Sturm attractor $\mathcal{A}_f$ of particular nonlinearities $f\in\mathcal{N}(x,u,u_x)$. The reason is purely on the ODE pendulum side~\eqref{C8:eq:2.1.2}. It is not an easy task to determine the equilibrium boundary order of $v\in\mathcal{E}_f$, as encoded in the Sturm permutation $\sigma_f$, for particular nonlinearities $f$. Indeed we address the whole class of dissipative nonlinearities $f=f(x,u,u_x)$, rather than any particular instance. \\

As a partial remedy we also consider the subclass $\mathcal{N}(u)$ of nonlinearities $f=f(u)$, for which the equilibrium ODE~\eqref{C8:eq:2.1.2} becomes $x$-autonomous Hamiltonian. See~\cite{FiedlerRochaWolfrum2012-SturmAttractors} and the survey~\cite{RochaFiedler14} for an explicit combinatorial characterization of this subclass of Sturm permutations $\sigma_f$. For a complete listing with $N=9$ equilibria see again Fig.~\ref{C8:fig:2.1.2}. The case $f\in\mathcal{N}(u,u_x^2)$ of nonlinearities $f(u,p) =f(u,-p)$ which are even in $p=u_x$ is only superficially more general.\ The equilibrium ODE~\eqref{C8:eq:2.1.2} is then only reversible in $x$, but not necessarily Hamiltonian. Still, reversibility implies integrability here. The ODE~\eqref{C8:eq:2.1.2}, but not the PDE~\eqref{C8:eq:2.1.1}, can be smoothly transformed to become Hamiltonian. The original Sturm permutation $\sigma_f$ becomes an involution, $\sigma_f=\sigma_f^{-1}$. Moreover $\sigma_f$ can be represented by a Hamiltonian nonlinearity $g\in\mathcal{N}(u)$ such that $\sigma_f=\sigma_g$. In particular, the Sturm global attractors $\mathcal{A}_f$ in the larger $x$-reversible class $\mathcal{N}(u,u_x^2)$ coincide with those of the $x$-reversible Hamiltonian class $\mathcal{N}(u)$. \\

Next, we turn to the case $f\in\mathcal{P}(x,u,u_x)$ of periodic boundary conditions $x\in S^1$. The nodal dropping property~\eqref{C8:eq:2.1.7} of zero numbers prevails. It has long been known, however, that any planar autonomous ODE vector field embeds into this larger class of PDEs~\cite{SandstedeFiedler92}. In particular, saddle-saddle heteroclinic and homoclinic trajectories may arise, and automatic transversality of stable and unstable manifolds fails. In the $S^1$-equivariant case of $x$-independent nonlinearities $f\in\mathcal{P}(u,u_x)$, on the other hand, automatic transversality does follow from mere (normal) hyperbolicity of equilibria and nonconstant rotating waves $u(t,x) = U\, (x-ct)$. In fact the global attractor then decomposes into equilibria $\mathcal{E}$, rotating waves $\mathcal{R}$, and their heteroclinic orbits $\mathcal{H}$ as
\begin{equation}\label{C8:eq:2.1.9}
\mathcal{A}_f = \mathcal{E}_f \cup \mathcal{R}_f \cup \mathcal{H}_f\,.
\end{equation}
In this periodic class, the PDE~\eqref{C8:eq:2.1.1} therefore defines a Morse--Smale system. In the $x$-reversible, i.e.\ $O(2)$-equivariant, class $f \in \mathcal{P}(u,u_x^2)$ an explicit Lyapunov function has been constructed; see~\cite{FiedlerGrottaRocha14}. Then all ``rotating'' waves are frozen, i.e.\ of wave speed $c=0$, and time-periodic solutions other than equilibria cannot arise. \\

In~\cite{FiedlerRochaWolfrum12}, we have developed an explicit algorithm to determine the heteroclinic connection graph in the general $S^1$-equivariant case $f_0\in\mathcal{P}(u,u_x)$. The basic procedure consisted in a homotopy $f_\tau$ from $f_0$ to $f_1\in\mathcal{P}(u,u_x^2)$ which preserved (normal) hyperbolicity of all equilibria and rotating waves, freezing the latter to wave speed $c=0$ at $\tau =1$. We were then able to derive the heteroclinic orbits in $\mathcal{P}$ from those in the reversible Neumann class $f_1\in \mathcal{N }(u,u_x^2)$ or, equivalently, the Hamiltonian Neumann class $f_1\in\mathcal{N}(u)$. \\

\begin{figure}
\centering \includegraphics[angle=0,scale=1]
{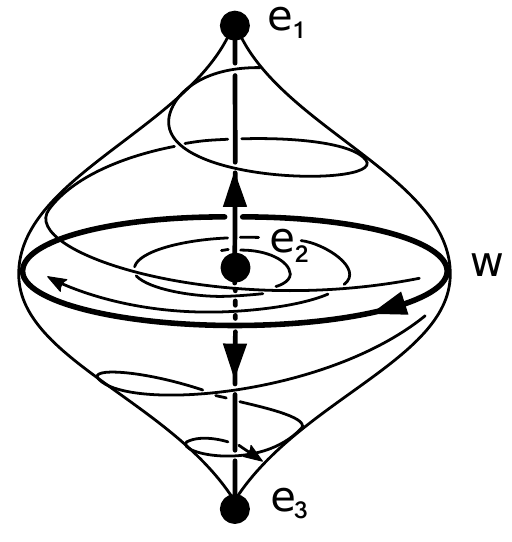}
\caption{\emph{A geometric representation of the Sturm spindle attractor with three equilibria and one rotating wave; see also~\cite{KrisztinWalther01}.}}
\label{C8:fig:2.1.3}
\end{figure}

As an explicit example we mention the \emph{Sturm spindle} attractor $\mathcal{A}$ of Fig.~\ref{C8:fig:2.1.3}, 
which consists of three equilibria $e_k$ of Morse indices $i(e_1)=i(e_3)=0$, $i(e_2)=3$ and one rotating wave $w$ with (strong) unstable dimension one. Surprisingly, the same global attractor $\mathcal{A}$ has also been encountered by Krisztin and Walther in the quite different setting of delay differential equations
\begin{equation}\label{C8:eq:2.1.10}
\dot{u}(t) = g\big(u(t),u(t-1)\big)\,,
\end{equation}
for positive monotone feedback $g$ in the delayed argument; see~\cite{KrisztinWalther01}. There is a nodal property similar to the zero number $z$ of~\eqref{C8:eq:2.1.8} at work, here.
The delay class is not exhaustively classified; it might turn out to be a subclass of the parabolic class, together with the closely related class of negative monotone delayed feedback. This striking analogy remains unexplored at present. \\

At the end of this section we return to the general Neumann class $f\in\mathcal{N}(x,u,u_x)$. We have a combinatorial description of this PDE class in terms of ODE shooting meanders, the Sturm permutations $\sigma_f$, and the explicit derivation of the heteroclinic connection graph $\mathcal{C}_f$. We did not yet succeed, however, to a priori characterize the class of all connection graphs $\mathcal{C}_f$ or, more ambitiously, the topological properties of the global Sturm attractors $\mathcal{A}_f$ themselves -- other than by mindless successive enumeration of all $\sigma_f$. However, some promising steps have been achieved. \\

We recall the gradient-like structure which implies a general decomposition $\mathcal{A}=\mathcal{E}\cup\mathcal{H}$ of the global attractor into (hyperbolic) equilibria $\mathcal{E}$ and their heteroclinic orbits $\mathcal{H}$. In other words,
\begin{equation}\label{C8:eq:2.1.11}
\mathcal{A}=\bigcup_{v\in\mathcal{E}}W^u(v)
\end{equation}
decomposes into the disjoint union of the unstable manifolds $W^u(v)$ of all equilibria $v$. Here $W^u(v)$ consists of all ancient trajectories $u(t,\cdot)$ in $\mathcal{A}$ which limit onto $v$, for $t\rightarrow-\infty$. In particular, $v \in W^u(v)$. Note that $\dim W^u(v)=i(v)$ is the Morse index of $v$, and $W^u(v)=\{ v\}$ for $i(v) = 0$. We call~\eqref{C8:eq:2.1.11} the \emph{dynamic complex} of $\mathcal{A}$. In general, the boundaries $\partial W^u =(\operatorname{clos}W^u)\smallsetminus W^u$ need not be manifolds. In fact,~\eqref{C8:eq:2.1.11} need not even be a cell complex, even though the $W^u(v)$ are embedded balls of dimensions $i(v)$. In the class of Sturm global attractors $\mathcal{A}= \mathcal{A}_f$, $f \in \mathcal{N}(x,u,u_x)$, however, the decomposition~\eqref{C8:eq:2.1.11} is a \emph{regular cell complex}. This means that the embedded $i(v)$-cells $W^u(v)$ possess embedded sphere boundaries $\partial W^u(v) \cong S^{i(v)-1}$, so that the attaching maps onto the boundary complexes are homeomorphisms, rather than just continuous maps. The proof requires a theorem of Schoenflies type for the embeddings of the closures $\operatorname{clos}W^u(v)$; see~\cite{Fiedler15}. For the resulting regular dynamic cell complex~\eqref{C8:eq:2.1.11} see~\cite{FiedlerRocha14}. \\

For a more concrete example, the case of a single closed Sturm 3-ball
\begin{equation}\label{C8:eq:2.1.12}
\mathcal{A}_f = \operatorname{clos}W^u(0) = W^u(0) \cup S^2\,.
\end{equation}
with trivial equilibrium $i(0)=3$ has been addressed in~\cite{FiedlerRocha14}. It has been shown how any finite regular cell decomposition of the Schoenflies boundary $\partial W^u(0) = S^2$ can be realized, topologically, as the dynamic complex
\begin{equation}\label{C8:eq:2.1.13}
S^2 = \bigcup_{v \in \mathcal{E} \smallsetminus \{0\}} W^u(v)
\end{equation}
of the remaining equilibria $v$ with Morse index $i(v) \leq 2$, in the Neumann class $f \in \mathcal{N}(x,u,u_x)$. \\

A more careful analysis of this example, however, also keeps track of the two-dimensional \emph{fast unstable manifold} $W^{uu}(0) \subset W^u(0)$ of those solutions $u(t, \cdot)$ in $W^u$ which do not decay to $v$ at the slowest possible exponential rate in $W^u(0)$, for $t \rightarrow -\infty$. By our Schoenflies results $\partial W^{uu}(0)=S^1$ decomposes $S^2$, as an ``equator'', into two (closed) hemispheres $S_\pm^2$. Each closed hemisphere, separately, is realizable as a planar Sturm attractor $\mathcal{A}_\pm = \mathcal{A}_{f_\pm}$ for some nonlinearities $f_\pm \in \mathcal{N}(x,u,u_x)$. To prove our Sturm realization~\eqref{C8:eq:2.1.13} in~\cite{FiedlerRocha14} we took one of the hemispheres to be a single closed 2-cell. The general question of admissible, i.e.\ Sturm realizable, hemisphere decompositions is unexpectedly delicate; see our work in progress~\cite{FiedlerRocha16a,FiedlerRocha16b,FiedlerRocha16c}. As a caveat, we just mention that the stable extreme equilibria $v_1, v_{27}$ in any dissipative Sturm realization of the 3d solid octahedron complex $\mathbb{O}$ of 27 equilibria must be chosen as adjacent corners of $\mathbb{O}$, rather than antipodally. For adjacent corners $v_1, v_{27}$, the two hemispheres $S_\pm^2$ are allowed to decompose $S^2$ into 1+7 or 2+6 octahedral faces, but not into 3+5 or 4+4 faces, by the fast unstable manifold $W^{uu}(0)$.\\

To complete our ambitious goal of a complete description of the Sturm regular dynamic complexes~\eqref{C8:eq:2.1.11}, an inductive process comes to mind. Suppose all Sturm global attractors $\mathcal{A}_f$, $f \in \mathcal{N}(x,u,u_x)$, have been understood for dimensions up to $n-1$. Two steps remain. In a first step we have to understand the sphere boundary of single Sturm balls $\mathcal{A}=\operatorname{clos}W^u(0)$, for $i(0)=n$, as gluings of closed hemisphere complexes $\mathcal{A}_\pm$ of dimension $n-1$, which are assumed to be known by induction hypothesis. In a second step we have to understand how these new $n$-cells attach to each other, and to lower-dimensional cells, to form a general Sturm global attractor of dimension $n$. Again, the respective hemisphere decompositions should play a crucial role here. At present, we have understood step~2 only up to the planar case $n=2$; see~\cite{FiedlerRocha2008-SturmAttractors2,FiedlerRocha2009-SturmAttractors1,FiedlerRocha2010-SturmAttractors3}. Step~1 for $n=3$ is imminent~\cite{FiedlerRocha16a,FiedlerRocha16b,FiedlerRocha16c}.

\subsection{Meanders}\label{subsec:meanders}

In the previous section, we introduced meander curves in the plane as shooting curves of parabolic PDEs~\eqref{C8:eq:2.1.1} in one space dimension; see in particular~\eqref{C8:eq:2.1.2},~\eqref{C8:eq:2.1.5} and Figs.~\ref{C8:fig:2.1.1},~\ref{C8:fig:2.1.2}. Quite independently from this ODE and PDE context, meanders also arise as trajectories of Cartesian billiards~\cite{FiedlerCastaneda2012-Meander}, and as representations of monomials in Temperley--Lieb algebras~\cite{DiFrancescoGolinelliGuitter1997-Meander}.
\begin{figure}
\centering \includegraphics[angle=-90,width=\textwidth]
{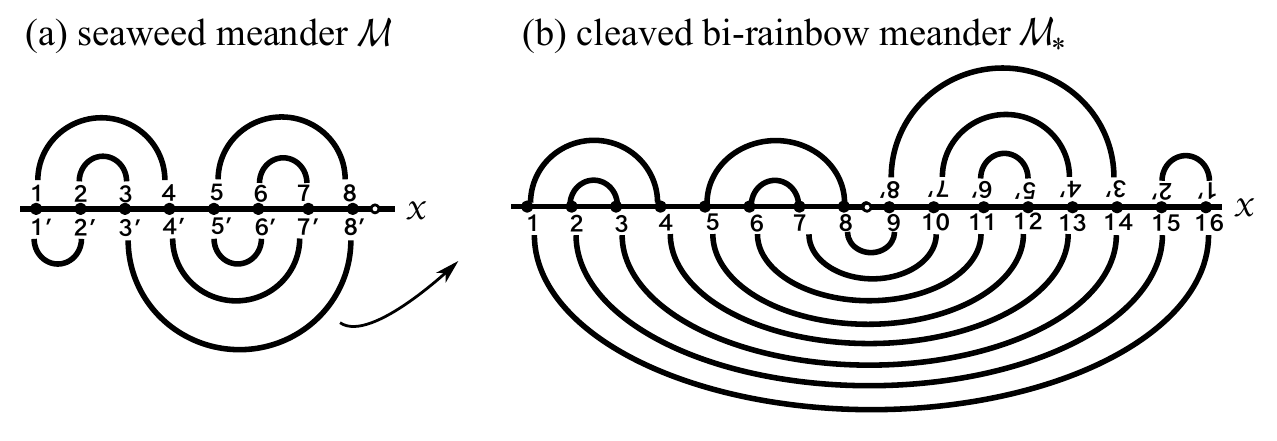}
\caption{\emph{(a) A seaweed meander $\mathcal{M}(2,2\,|\,1,3)$ with $N=8$ intersections.
(b) The associated cleaved bi-rainbow meander $\mathcal{M}_*(2,2,3,1\,|\,8)$ with $2N=16$ intersections.
Note how upper arches $(i,j)$ of $\mathcal{M}$ have been preserved.
Lower arches $(i,j)$ of $\mathcal{M}$ have been converted to upper arches $(2N+1-i,\, 2N+1-j)$ of $\mathcal{M}_*$; see~\cite{FiedlerCastaneda2012-Meander}.}}
\label{C8:fig:2.2.1}
\end{figure}
For simplicity we now represent meander curves, or equivalently meander permutations, by $\left[N/2\right]$ standardized upper and lower half-circular arches joining at vertices $1,\dots,N$ along the horizontal axis. Also we replace the open meanders from south-west to north-east, of Section~\ref{subsec:sturmglobalattractors}, by their closed Jordan curve counterparts. See Fig.~\ref{C8:fig:2.1.1}(c) for one simple closing procedure. \\

More generally, we now regard closed meanders $\mathcal{M}$ as the pattern created by one or several disjoint closed smooth Jordan curves in the plane as they intersect a given horizontal line transversely. Let $c(\mathcal{M}) \geq 1$ count the connected components of $\mathcal{M}$. In view of Sturm global attractors, of course, the case $c(\mathcal{M}) =1$ of connected meanders generated by a single closed Jordan curve is of primary interest. \\

Several new results have been obtained which are based on rainbows. Here a \emph{proper upper rainbow} consists of a number $\alpha$ of concentric and adjacently nested upper arches in any meander. It is understood that $\alpha$ is chosen maximal, in all cases. An \emph{upper rainbow} consists of one or several neighboring proper upper rainbows. Lower rainbows are defined analogously. Note that upper and lower rainbows may themselves be nested, giving rise to a forest of their sizes $\alpha$. A closed meander without rainbow nesting is called a \emph{seaweed meander}. Symbolically a seaweed meander can be described as
\begin{equation}\label{C8:eq:2.2.1}
\mathcal{M}(\alpha_1,\dots,\alpha_{m_+}\,|\, \beta_1,\dots,\beta_{m_-})\,,
\end{equation}
where $\alpha_1,\dots, \alpha_{m_+}$ are the sizes of the neighboring proper upper rainbows, and $\beta_1,\dots, \beta_{m_-}$ enumerate the proper lower rainbows, left to right. See Fig.~\ref{C8:fig:2.2.1}(a) for an example. Trivially
\begin{equation}\label{C8:eq:2.2.2}
\alpha_1+\dots+\alpha_{m_+}=\beta_1+\dots+\beta_{m_-}=N/2
\end{equation}
for any closed meander of $N$ vertices. The special seaweed case $m_- =1$ of a single lower rainbow of size $\beta = N/2$ is called a \emph{bi-rainbow meander}.
\begin{figure}
\centering \includegraphics[angle=-90,width=0.8\textwidth]
{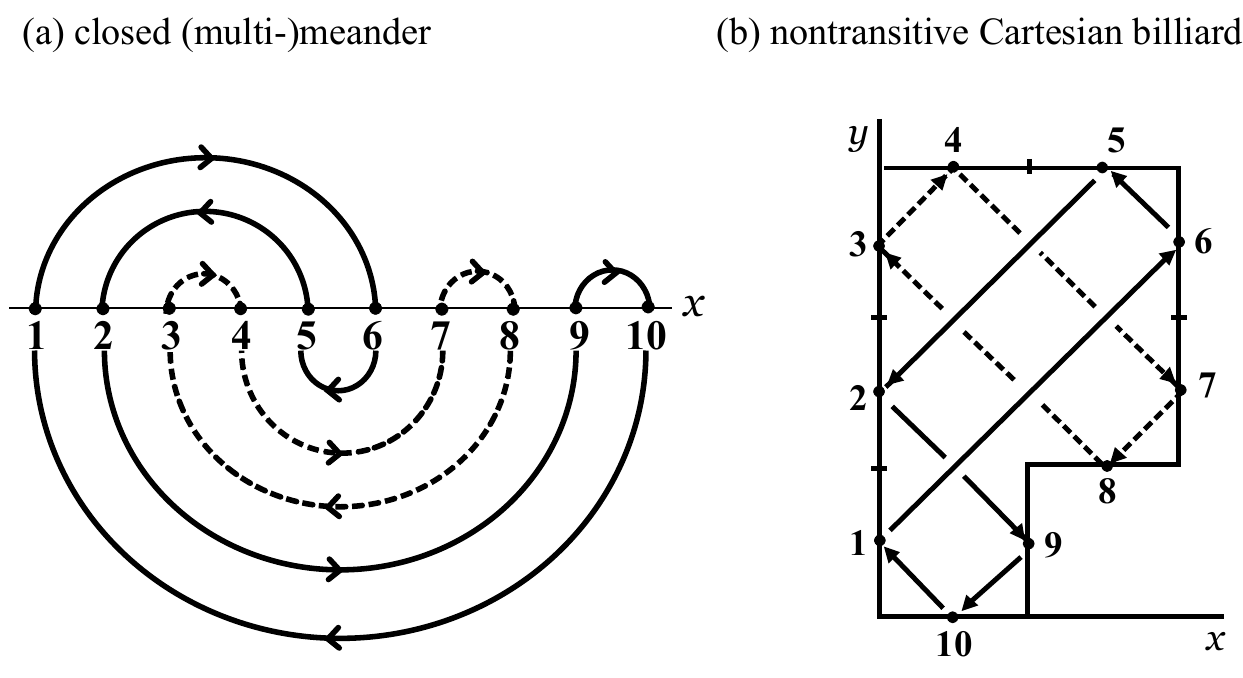}
\caption{\emph{(a) A bi-rainbow meander with two connected components, one distinguished by dashed arches.
(b) A corresponding nontransitive plane Cartesian billiard. One closed flight path is distinguished by dashed lines; see~\cite{FiedlerCastaneda2012-Meander}.}}
\label{C8:fig:2.2.2}
\end{figure}
The case where the upper rainbows may be arbitrarily nested, but there is only one single lower rainbow $\beta = N/2$, is called \emph{lower rainbow meander}; see Fig.~\ref{C8:fig:2.2.1}(b). For example any closed meander with $N$ vertices can be ``opened'' to a closed lower rainbow meander with $2N$ vertices and the same number of connected components. Indeed we just duplicate the original vertices $j \mapsto j,j'$ and open up the horizontal axis to accommodate the vertices in the order $1,\dots, N,N',\dots, 1'$; see again Fig.~\ref{C8:fig:2.2.1}. Recall how only connected meanders relate to an interpretation as ODE shooting curves for PDE Sturm global attractors. Lower rainbow meanders of the above general type were related to Cartesian billiards by Fiedler and Casta\~{n}eda; see~\cite{FiedlerCastaneda2012-Meander} and Fig.~\ref{C8:fig:2.2.2}. \\

In fact, there exists a plane with only horizontal and vertical billiard boundaries, with vertices all on the Cartesian integer grid and billiard slopes of $\pm 45^\circ$ between half-integer boundary points, such that the connectivities of paths in the billiard and in the meander coincide. As an almost immediate consequence one obtains the following list for the number $c(\alpha_1,\dots,\alpha_m)$ of connected components of bi-rainbows $\mathcal{M}(\alpha_1,\dots, \alpha_n\,|\,\beta)$ in terms of greatest common divisors $\gcd$:
\begin{equation}\label{C8:eq:2.2.3}
\begin{aligned}
c(\alpha_1) &= \alpha_1\,,\\
c(\alpha_1,\alpha_2) &= \gcd (\alpha_1, \alpha_2)\,,\\
c(\alpha_1, \alpha_2, \alpha_3) &= \gcd (\alpha_1 +\alpha_2, \alpha_2 + \alpha_3)\,.
\end{aligned}
\end{equation}
The dissertation of Karnauhova~\cite{Karnauhova2016-Dissertation} establishes how such formulas cannot be extended to $m\geq4$. In fact $c$ cannot be written as the $\gcd$ of any two polynomials $f_1,f_2$ in $\alpha_1,\dots,\alpha_m$ with integer coefficients. In~\cite{KarnauhovaLiebscher2015-Meanders} this obstacle is overcome, replacing the algorithm ``$\gcd$'' by an explicit recursion of related logarithmic complexity in terms of the data $\alpha_1,\dots,\alpha_m$. \\

\begin{figure}
\centering \includegraphics[angle=0,width=\textwidth]
{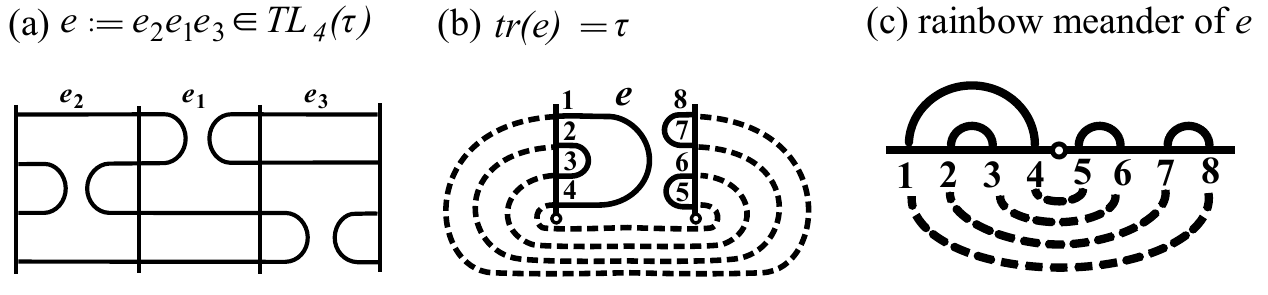}
\caption{\emph{(a) Diagram representation of the monomial $e=e_2e_1e_3 \in TL_4(\tau)$.
(b) Homotopy equivalent diagram of $e$ with dashed exterior strands matching right and left ends of the same strand index.
Note the number $k=1$ of exterior connected components of the closed diagram.
(c) Joining the lower ends ``$o$'' of the vertical strand boundaries of $e$, and opening the boundaries to become horizontal we obtain an equivalent lower rainbow meander with $8$ intersections and the same number of connected components as in (b); see~\cite{FiedlerCastaneda2012-Meander}.}}
\label{C8:fig:2.2.4}
\end{figure}

Encouraged by suggestions of Matthias Staudacher (projects C4 ``Structure of Quantum Field Theory: Hopf Algebras versus Integrability'' and C5 ``AdS/CFT Correspondence: Integrable Structures and Observables'') and earlier hints by Eberhard Zeidler, we also began to explore relations between meanders and Temperley-Lieb Hopf algebras $TL_N(\tau)$; see~\cite{FiedlerCastaneda2012-Meander}. Algebraically these are given by $N$ multiplicative generators $e_0=1,e_2,\dots,e_{N-1}$ with the three relations
\begin{equation}\label{C8:eq:2.2.4}
\begin{aligned}
e_i^2 &=\tau e_i\,,\\
e_i e_j &= e_j e_i\,,\\
e_i e_{i \pm 1} e_i &= e_i\,,
\end{aligned}
\end{equation}
for all appropriate $i,j >0$ with $|i-j|\, >1$.
\begin{figure}
\centering \includegraphics[angle=0,width=\textwidth]
{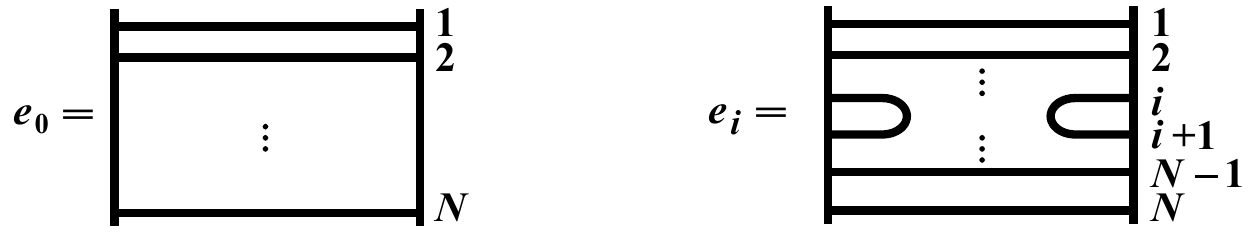}
\caption{\emph{Diagram representation of Temperley--Lieb generators $e_0=1$ and $e_i$, for $i=1,\dots,N-1$; see~\cite{FiedlerCastaneda2012-Meander}.}}
\label{C8:fig:2.2.3}
\end{figure}
See Fig.~\ref{C8:fig:2.2.3} for a representation of the Temperley--Lieb generators $e_i \in TL_N(\tau)$ by $N$-strand diagrams. Multiplication of generators is realized by horizontal concatenation of the generator diagrams, up to strand isotopy. Isolas effect multiplication by $\tau$, each. As illustrated in Fig.~\ref{C8:fig:2.2.4} for $e$:$=e_2e_1e_3 \in TL_4(\tau)$, Temperley--Lieb monomials $e$ correspond to lower rainbow meanders $\mathcal{M}[e]$. The number $c$ of connected components of $\mathcal{M}[e]$ is then related to the \emph{Markov trace} $\operatorname{tr}_F$ of $e$:
\begin{equation}\label{C8:eq:2.2.5}
\operatorname{tr}_F(e) = \tau^{c(\mathcal{M}[e])}\,,
\end{equation}
see~\cite{DiFrancescoGolinelliGuitter1997-Meander}. For further comments and a relation to the classical Young-Baxter equations see also the dissertation by Karnauhova~\cite{Karnauhova2016-Dissertation} and the references there. 
\\

The case~\eqref{C8:eq:2.2.1},~\eqref{C8:eq:2.2.2} of seaweed meanders $\mathcal{M}(\alpha_1,\dots, \alpha_{m_+}\,|\, \beta_1,\dots,\beta_{m_-})$ is related to the index of seaweed algebras, as studied by Dergachev and Kirillov~\cite{DK00}. The seaweed algebra of type $(\alpha_1,\dots,\alpha_{m_+}\,|\, \beta_1,\dots,\beta_{m_-})$ is the $\operatorname{gl}(N/2)$ subalgebra of matrices which preserve the vector spaces $V_j$ spanned by unit vectors $e_1,\dots,e_{\alpha_1+\ldots+\alpha_j}$ and $W_j$ spanned by $e_{\beta_1+\ldots+\beta_j+1},\dots,e_{N/2}$. In fact the seaweed meander corresponds to the Cartesian billiard on the domain of admissible nonzero entries of the seaweed matrices. The index of the algebra is the minimal kernel dimension of the map $A \mapsto f[A,\cdot]$, for any linear functional $f$ in the dual algebra. By~\cite{DK00} this index coincides with the number $c(\mathcal{M})$ of connected components of the meander $\mathcal{M}$. \\

In her dissertation, Karnauhova considers the maximally disconnected case $\alpha_j=\beta_j$ of seaweed algebras with maximal index $N/2$, generated by trivial $2$-cycles of $\sigma$, for the closed meander $\mathcal{M}$; see~\cite{Karnauhova2016-Dissertation}. She then shifts the lower $\beta$-rainbows one step to the right. Ratcheting into the upper $\alpha$-rainbows, she obtains a dissipative Morse meander $\sigma_f$, and hence a Sturm global attractor $\mathcal{A}_f$. The Chafee--Infante attractors \cite{Henry85}, for example, which arise for symmetric cubic nonlinearities $f$ and have been the first Sturm global attractors studied historically, correspond to $m_\pm=1$ and $\alpha=\beta$. They are the Sturm global attractors of maximal dimension $\alpha =N/2$ for any given (odd) number $N+1$ of equilibria. The same construction in fact works for general, but identical, upper and lower arch configurations. Although the resulting global attractors are all pitchforkable, and thus have been known ``in principle'' by Fusco and Rocha~\cite{FuscoRocha}, the new approach is elegant and geometric: the Sturm global attractors can be obtained explicitly by successive double cone suspensions and elementary gluing concatenations. \\

In~\cite{Fiedler14}, Fiedler embarked on a different approach to explore the relation between the chaotic dynamics of ($x$-periodically forced) equilibrium equations~\eqref{C8:eq:2.1.2} with nonlinearity $f$ and the Sturm permutations $\sigma_f$ of the resulting shooting meanders. Homoclinic tangles with their associated shift dynamics are taken as a paradigm. In this direction, meander permutations have been defined for standard hyperbolic linear $SL_2(\mathbb{Z})$ torus automorphisms. The meander permutations compare the orderings of homoclinic points along the one-dimensional stable and unstable manifolds of their limits, respectively. As such, they are new topological invariants -- with intriguing relations to continued fraction expansions and quadratic number fields.

\subsection{Stability, instability, and obstacles}\label{subsec:statsol}

Consider the higher-dimensional analogue of~\eqref{C8:eq:2.1.1}, that is, for semilinear equations
$$u_t=\Delta u+f(u)$$
on a domain $\Omega\subseteq{\mathbb R}^N$ with $N\ge1$. More generally consider systems $u\colon\Omega\to{\mathbb R}^n$
\begin{equation}\label{C8:eq:semilinear}
u_t=Lu+f(u)
\end{equation}
with an elliptic differential operator $L$ and associated Dirichlet, Neumann, or mixed boundary values. Then the naturally associated function spaces $W_\Gamma^{1,2}(\Omega)$ and $L^2(\Omega)$ are not continuously embedded into $C(\overline\Omega)$. Therefore, differentiability of the corresponding superposition operator $F(u)(x)=f(u(x))$ fails in $L^2(\Omega)$ (except in the linear case). This imposes an unrealistically restrictive growth condition for $f$ in case $F\colon W_\Gamma^{1,2}(\Omega)\to L^2(\Omega)$. Moreover, no analogue of the zero number~\eqref{C8:eq:2.1.7} with a property like~\eqref{C8:eq:2.1.8} is available. Nevertheless, we are interested in the long-term behavior of the equation in these two function spaces. \\

A particularly natural question is whether a ``formally'' linearized stability of some stationary solution $u_0$ implies asymptotic stability in these spaces. In other words, suppose the spectrum of the ``formal'' linearization $L+f'(u_0)$ lies in the complex left half-plane. Do all solutions starting close to $u_0$ in the topology of the space $W_\Gamma^{1,2}(\Omega)$ or $L^2(\Omega)$ indeed remain close to $u_0$ and eventually tend to $u_0$ in the corresponding topology? Let $A$ denote the operator naturally associated to $-L$ and the boundary conditions. Suppose $A$ solves Kato's square root problem, that is, if $A^{1/2}\colon W_\Gamma^{1,2}(\Omega)\to L^2(\Omega)$ is an isomorphism. Then $W_\Gamma^{1,2}(\Omega)=D(A^{1/2})$. The asymptotic stability in $W_\Gamma^{1,2}(\Omega)$ is known for $F\colon W_\Gamma^{1,2}(\Omega)\to L^2(\Omega)$; see e.g.~\cite{Henry}. Gurevich and V\"{a}th~\cite{GurevichVaeth} have relaxed this hypothesis to a subcritical growth condition on $f$. Moreover, they obtained an asymptotic stability result in $L^2(\Omega)$. The method of proof was to consider extrapolated spaces and extrapolated semigroups which were obtained by two different procedures. As a side result, they have shown that these extrapolation procedures coincide if and only if $A$ solves Kato's square root problem. \\

One particular motivation to consider the topology of the spaces $W_\Gamma^{1/2}(\Omega)$ and $L^2(\Omega)$ is that the following question so far can be addressed only in these topologies. Suppose that the problem~\eqref{C8:eq:semilinear} is endowed with some obstacle, i.e., the values of $u$ on some parts of $\Omega$ or of its boundary are forced to stay within a certain cone by some constraint. One of our surprising results is that the principal of linearize asymptotic stability may then fail in $W_\Gamma^{1,2}(\Omega)$. The proof is based on an instability criterion for variational inequalities and a rather general formula relating the fixed point index of the flow of the obstacle problem with the topological degree of a certain associated operator~\cite{KimVaeth}. For specific examples see~\cite{KimVaeth,VaethProc}.
\begin{figure}
\centering \includegraphics[angle=0,width=\textwidth]
{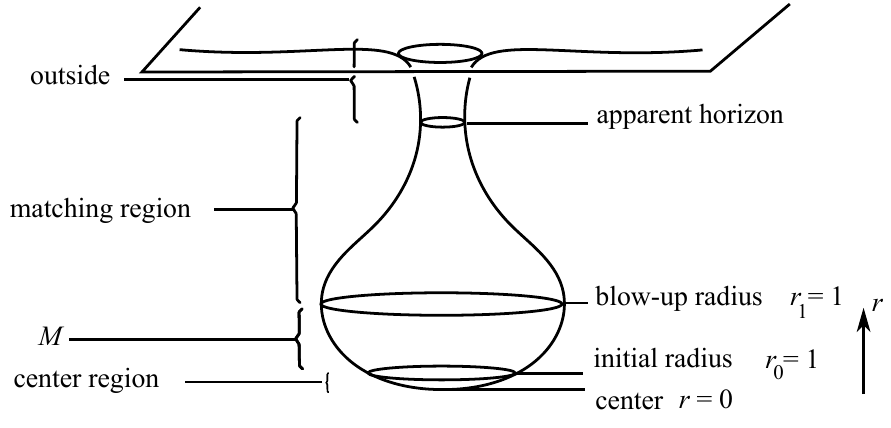}
\caption{\emph{A $2d$-caricature of black hole initial data for the Einstein equations. The vertical variable is the radius $r$. Each horizontal section $\{ r=\operatorname{constant}\}$ in this figure is a circle, which represents a $2$-sphere $S^2$. The region outside the apparent horizon (and above it in this picture) was treated in the literature. In~\cite{FiedlerHellSmith15} we construct the region $M$ between an initial radius $r_0>0$ and a blow-up radius $r_1 =1$ whose section $\{r_1\} \times S^2$ is a critical point of the area functional.\ The matching and center regions remain to be constructed; see~\cite{FiedlerCastaneda2012-Meander}.}}
\label{C8:fig:2.3.1}
\end{figure}
\begin{figure}
\centering \includegraphics[angle=-90,width=\textwidth]
{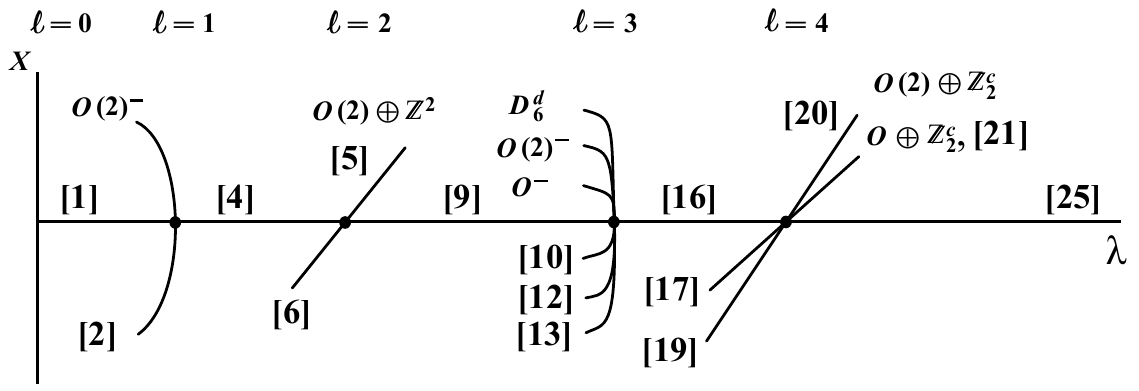}
\caption{\emph{Local transcritical and pitchfork bifurcation branches $(\lambda_K(s), v_K(s))$ from $\lambda_K(0) = \lambda_\ell=\ell(\ell+1)$, $v_K(0)=0$ with maximal isotropies $K$ as indicated.
Numbers in brackets indicate co-dimensions of the strong unstable manifolds of their normally hyperbolic group orbits; see~\cite{FiedlerHellSmith15}.}}
\label{C8:fig:2.3.2}
\end{figure}
\begin{figure}
\centering \includegraphics[angle=0,width=\textwidth]
{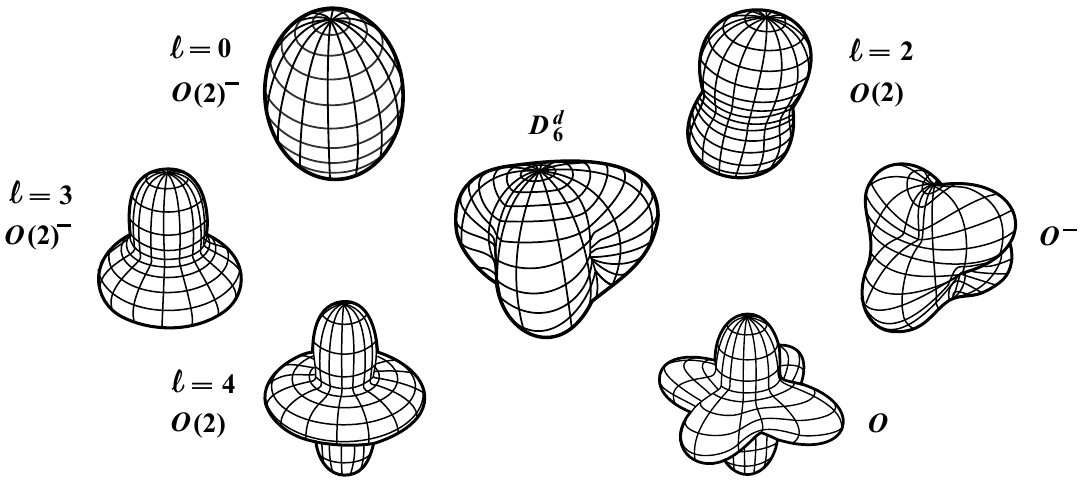}
\caption{\emph{The graphs of the bifurcating symmetry breaking equilibria $v_K$ over the sphere $S^2$, with isotropies $K$ as indicated; see~\cite{FiedlerHellSmith15}.}}
\label{C8:fig:2.3.3}
\end{figure}
%

\subsection{Symmetric and asymmetric black hole data}\label{subsec:blackholedata}

The Cauchy initial-value problem for the relativistic Einstein equations requires initial data which satisfy the Einstein constraint.
For black holes these constraints take the form of an equation for the metric $g$ with prescribed scalar curvature $R$ on an asymptotically flat Riemannian $3$-manifold $M$; see Fig.~\ref{C8:fig:2.3.1} for an illustration and~\cite{Ringstroem2009-CauchyGeneralRelativity} for a general background. For simplicity we consider a collar region $M= [r_0, r_1] \times S^2$ with standard $2$-sphere cross section. In a suitable coordinate foliation, the metric can be expressed as
\begin{equation}\label{C8:eq:2.3.1}
g=u^2dr^2 + r^2\omega
\end{equation}
with the standard metric $\omega$ on $S^2$. For example consider prescribed scalar curvatures $R$ such that $\lambda$:= $r^2 R -2$ is a constant bifurcation parameter. Then $u=u(r,x)$ can be rescaled to
\begin{equation}\label{C8:eq:2.3.2}
u(r,x) = \big( \tfrac{\lambda}{2} (\tfrac{1}{r}-1) \big)^{-1/2}
\big( v \big(-\tfrac{1}{2} \log (1-r),x\big) +1 \big)\,,
\end{equation}
and the Einstein constraints on the radial metric coefficient $u$ take the quasilinear parabolic form
\begin{equation}\label{C8:eq:2.3.3}
v_t = (1+v)^2 \big( \Delta v + \lambda f(v)\big)\,.
\end{equation}
Here $v = v(t,x)$, $x \in S^2$, $\Delta$ is the Laplace-Beltrami operator on $S^2$, and $f(v) = v- \tfrac{1}{2}v^2/(1+v)$ defines a cubic nonlinearity. \\

In~\cite{FiedlerHellSmith15}, we obtain spherically anisotropic solutions $v$ for isotropically prescribed scalar curvature $R$. Equilibria $v(b,x)=v_*(x)\not\equiv\operatorname{const}$ provide self-similar anisotropic metrics $g$. Heteroclinic trajectories $u(t,x)\rightarrow v_\pm(x)$, for $t \rightarrow\pm\infty$, are asymptotically self-similar for $r \rightarrow 0$ and $r \rightarrow 1$. The metric $g$ remains smooth, in fact, in the latter limit. We obtain these solutions by an $O(3)$ symmetry breaking bifurcation analysis of the scaled quasilinear PDE in a center manifold at the trivial isotropic equilibrium $v_*\equiv0$. \\

Symmetry breaking bifurcations occur at $\lambda_\ell = \ell(\ell+1)$, for $l=0,1,2,\dots$, in the spherical harmonics kernels $V_\ell$ of dimension $2\ell+1$. See Fig.~\ref{C8:fig:2.3.2} for the bifurcating branches of nonconstant equilibria $v$, and Fig.~\ref{C8:fig:2.3.3} for their remaining reduced, nonspherical spatial isotropies. Heteroclinic trajectories follow, between $v_*$ and $0$, by standard bifurcation theory in the one-dimensional isotropy subspaces of each kernel $V_\ell$. For heteroclinic trajectories between equilibria of different isotropies; see also~\cite{ChossatLauterbach}. \\

The results above are of a local nature, near the spherically isotropic trivial equilibrium $v_* \equiv 0$. The fully isotropic case $\ell=0$ was treated in~\cite{smith2007blow} previously, by more direct methods. In practice specific computations, particularly of the unstable dimensions and the bifurcation directions, are limited to low $\ell$ and have been carried out for $\ell\leq4$, only. Axisymmetric solutions $v=v(t,\vartheta)$, however, which depend only on the azimuthal angle $\vartheta$ of $x \in S^2$, are known to bifurcate for all $\ell=0,1,2,\dots$. Clearly~\eqref{C8:eq:2.3.3} becomes ``spatially'' one-dimensional,
\begin{equation}\label{C8:eq:2.3.4}
v_t=(1+v)^2 \big(v_{\vartheta \vartheta} + \cot \vartheta\, v_\vartheta + \lambda f(v) \big)\,,
\end{equation}
with Neumann boundary on $0 < \vartheta < \pi$. The dissertation of Lappicy explores these versions of quasilinear Sturm parabolic equations with a singular spatial coefficients, by the global methods of Section~\ref{subsec:sturmglobalattractors}; see~\cite{Lappicy2016-Dissertation}. In particular, this yields a plethora of axisymmetric, but not fully spherically isotropic, equilibrium solutions and their heteroclinic trajectories as parts of black hole initial data which satisfy the Einstein constraint.

\section{Blow-up and Bianchi cosmology}

\numberwithin{equation}{section}
\numberwithin{figure}{section}

\subsection{Dynamics at infinity and blow-up in complex time}

In the previous Section~\ref{sec:globalattractors}, a class of dissipative equations~\eqref{C8:eq:2.1.1} was studied, where longtime existence and boundedness of the solutions was guaranteed. In the present section we focus on equations which develop singularities in finite or infinite time. In other words, some trajectories escape to infinity. In order to study the dynamics at infinity, one possible approach uses \emph{Poincar\'e ``compactification''}: the solution space is projected centrally onto a tangent upper hemisphere. More precisely the hemisphere at its north pole is tangent to the solution space at its origin. Arbitrarily far points in the solution space are thereby mapped to the equator of the hemisphere, also called the \emph{sphere at infinity}. For details see~\cite{helldiss,hell13,perko}. We keep the name Poincar\'e ``compactification'' which originates from ODE theory, where the solution space is $\mathbb{R}^n$ and the compactified space is a compact hemisphere. In the PDE context, the solution space is infinite dimensional and its ``compactification'' is a bounded hemisphere, which is not compact. \\

As an elementary but instructive example, consider the linear equation: 
\begin{equation}\label{C8:linearreacdiff}
u_t=u_{xx}+bu\,,\\
\end{equation} 
for $0<x<\pi$, with Neumann boundary conditions; see~\cite{helldiss}. The equation~\eqref{C8:linearreacdiff} is not dissipative for positive $b$: some trajectories escape to infinity, in infinite time. The class of equations at the boundary between dissipativity and finite time blow-up is called \emph{slowly nondissipative} in~\cite{BGdiss}: they admit at least one immortal trajectory escaping to infinity in forward time. The linear example~\eqref{C8:linearreacdiff} admits only one bounded equilibrium, located at the origin, and generates infinitely many equilibria at infinity. Those are located at the intersections of the one-dimensional eigenspaces, associated to the eigenfunctions of the Laplacian, with the sphere at infinity. Let us denote the equilibria at infinity by $\pm\Phi_m$, $m\in\mathbb{N}$. For all $b\neq n^2$ which are not square integers, all equilibria, both bounded and at infinity, are hyperbolic. Furthermore the dynamics within the sphere at infinity has a Chafee--Infante structure; see Section~\ref{subsec:sturmglobalattractors}. More precisely the two equilibria $\pm\Phi_0$ at infinity are stable. The equilibria $\pm\Phi_{m}$ at infinity each possess an $m$-dimensional unstable manifold that consists of heteroclinic orbits to the equilibria $\pm\Phi_{j}$, $j\in\{0,\dots,m-1\}$. The unstable manifold of the origin is of finite dimension $[\sqrt{b}]+1$, where $[\sqrt{b}]$ denotes the integer part of $\sqrt{b}$. This unstable manifold consists of heteroclinic orbits to the equilibria $\pm \Phi_i$, $i\in\{0,\dots,[\sqrt{b}]\}$. The stable manifold of the equilibrium at the origin is infinite dimensional and consists of heteroclinic orbits from the equilibria $\pm\Phi_m$ at infinity, $m>[\sqrt{b}]$, to the origin. This describes the global unbounded attractor for~\eqref{C8:linearreacdiff}. \\

When~\eqref{C8:linearreacdiff} is augmented by a smooth bounded nonlinearity $g=g(u)$, 
\begin{equation}\label{C8:linreacdiff+bd}
u_t=u_{xx}+bu+g\,,\\
\end{equation}
the dynamics at infinity remains essentially unchanged. In her dissertation~\cite{BGdiss}, Ben-Gal adapts the techniques of Section~\ref{subsec:sturmglobalattractors} to study the global heteroclinic structure of~\eqref{C8:linreacdiff+bd}. In particular, she proves the existence of a finite-dimensional unbounded inertial manifold containing the so called \emph{unbounded} (or \emph{noncompact}) \emph{attractor}. This object consists of all bounded equilibria (generically hyperbolic), their finite dimensional unstable manifolds, the finitely many equilibria at infinity (also generically hyperbolic) to which those unstable manifolds connect, and the finite dimensional unstable manifolds of those equilibria at infinity. Theorems similar to the dissipative case reveal the heteroclinic structure of the unbounded global attractor of slowly nondissipative equations. Pimentel and Rocha generalized some of these results to bounded nonlinearities $g(x,u,u_x)$; see~\cite{pim14, pim15}. \\

We consider blow-up in finite time for~\eqref{C8:eq:2.1.1}, $f=f(u)$, next. For smooth right hand side with polynomial growth, time can be rescaled such that trajectories do not run into equilibria at infinity in finite time. The resulting semiflow on the closed Poincar\'e hemisphere is nonsingular. Conley index theory then provides heteroclinic orbits between equilibria, or general isolated invariant sets. Isolated invariant sets are contained in isolating neighborhoods. An isolating neighborhood is a closed neighborhood whose boundary points eventually leave the neighborhood in forward or backward time direction. For details see~\cite{helldiss} and the references there. \\

It turns out that even for some polynomial ordinary differential equations, equilibria at infinity might not be isolated invariant sets. 
\begin{figure}
\begin{center}
\includegraphics[width=0.3\textwidth]{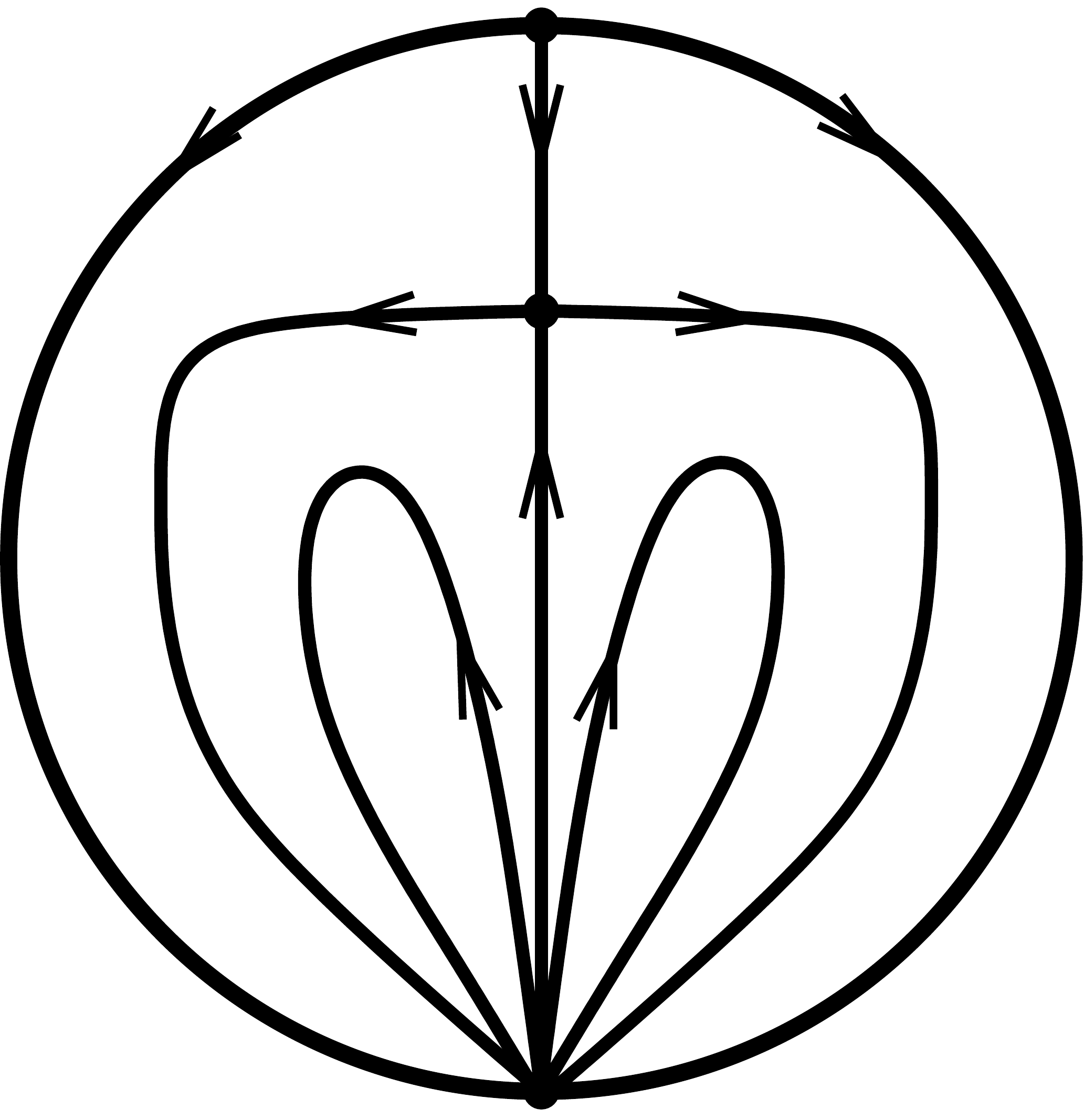}
\hspace{2cm}
\includegraphics[width=0.28\textwidth]{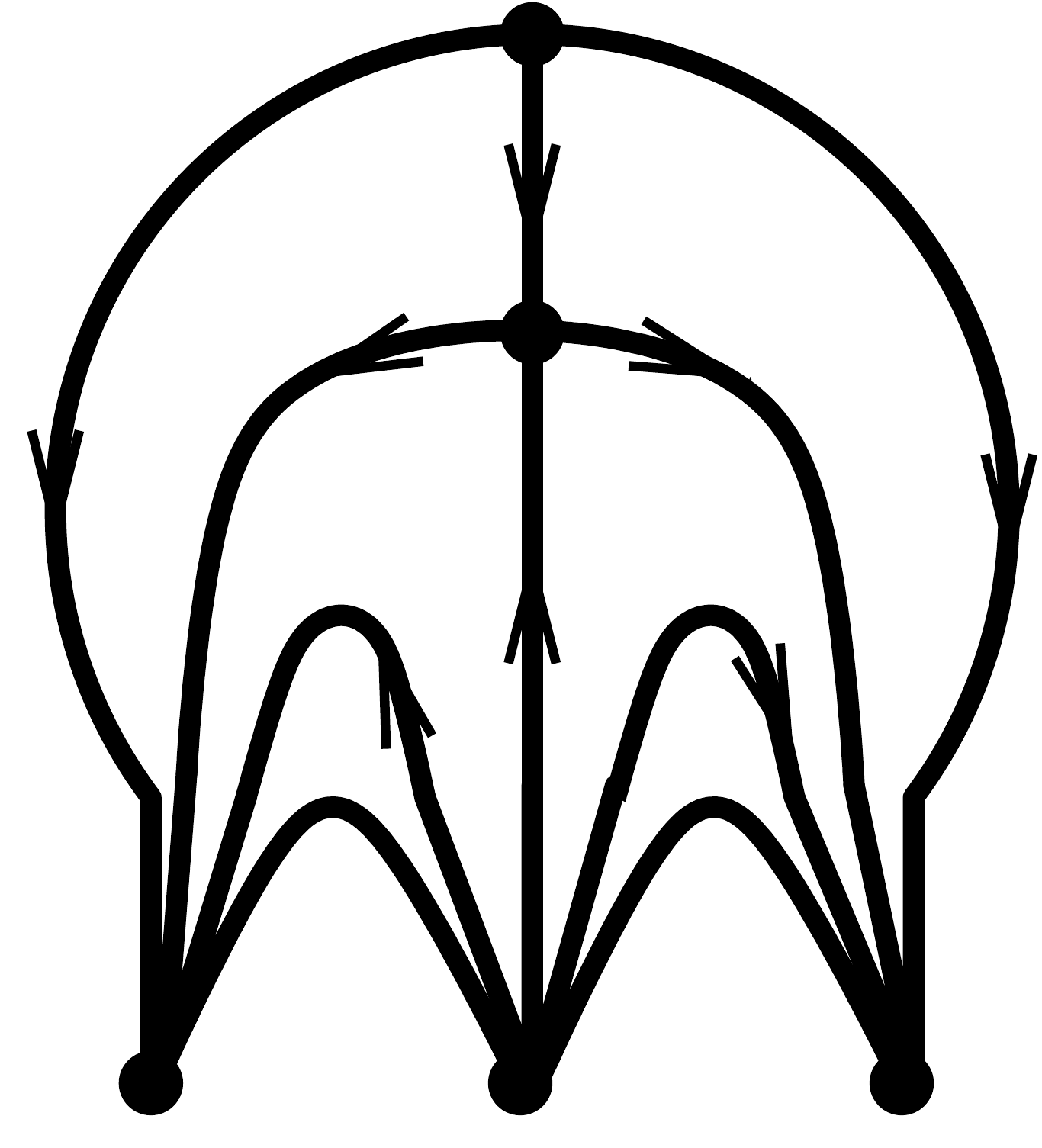} 
\end{center}
\caption{\emph{Left: the lower equilibrium at infinity is not an isolated invariant set, but is of isolated invariant complement. Right: the ersatz infinity for that equilibrium.}}
\label{homflower}
\end{figure}
Fig.~\ref{homflower}, left hand side, illustrates such an equilibrium at infinity. This equilibrium is isolated, invariant, and stable within the sphere at infinity. Towards the interior of the Poincar\'e hemisphere, however, it shows two petals of nested homoclinic loops. Therefore it does not possess small isolating neighborhoods. The complement of any small neighborhood of the degenerate equilibrium, however, is an isolating neighborhood. An invariant set with this property is called of \emph{ isolated invariant complement}. Although the complementary invariant set isolated by such a neighborhood is possibly very complex, its Conley index is well-defined. Poincar\'e duality then defines a meaningful Conley index for the degenerate equilibrium at infinity. An invariant set at infinity attracts some finite initial conditions, if and only if it is not a repeller. This can be characterized in terms of Conley index. Therefore this tool allows us to determine which part of the dynamics in the sphere at infinity is reachable by grow-up or blow-up solutions. \\

Via a standard construction illustrated in the right hand side of Fig.~\ref{homflower}, and described in detail in~\cite{helldiss,hell13}, it is possible to replace the invariant set of isolated invariant complement by a benign \emph{ersatz} with a well-defined Conley index. This ersatz can be used for the analysis of the global heteroclinic structure via classical Conley index theory. \\
 
Based on the idea of self-similar scaling and the Poincar\'e sphere at infinity, Fiedler has contributed to an unexpected application in a certain renormalization approach to the experimentally observed formation of spin liquids in carbon nanotubes; see~\cite{hsufied}. Renormalization has been suggested to overcome the limitations of mean-field descriptions. In its most simple form, homogeneously quadratic ODE systems
\begin{equation}
\frac{dx}{dt}=-\nabla_{x}V
\end{equation}
were addressed. The variational structure features a cubic functional
\begin{equation}
V(x)=-\frac13 \sum_{i=1}^{N} x_i (x^T A_i x)
\end{equation}
and real symmetric $N\times N$-matrices $A_i$. As a result we obtain self-similar and asymptotically self-similar blow-up solutions, as well as their blow-up rates. These results provide a rational foundation for the renormalization group approach to large, strongly correlated electron systems. \\

\begin{figure}
\centering \includegraphics[angle=0,scale=1.1]{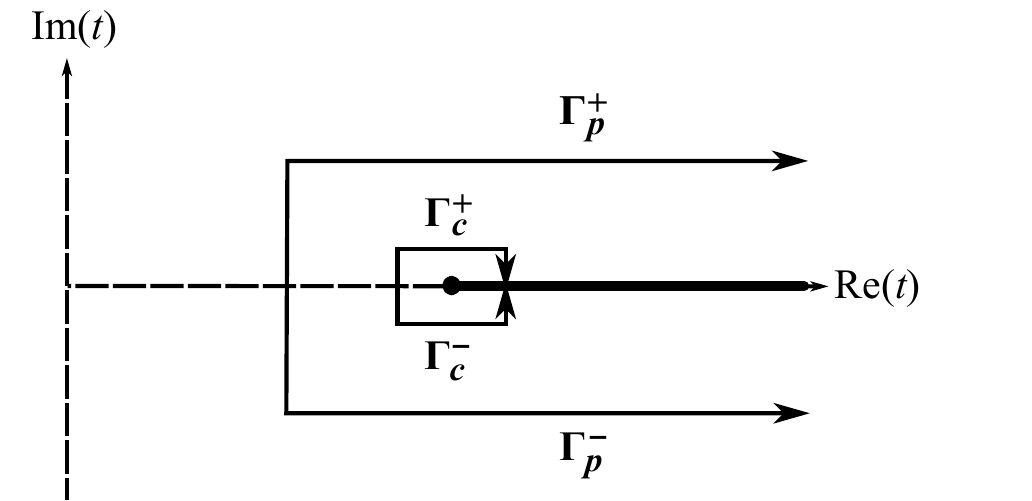}
\caption{\emph{The complex time domain in which we consider the continuation $\Gamma^\pm_p$, $\Gamma^\pm_c$ of blow-up solutions. The dot indicates the blow-up time $ T$. After $t=T$ we introduce a cut along the positive real axis, drawn as a bold line.}}\label{fig:complex_time_domain}
\end{figure}

For fast nonlinear diffusion in several space dimensions the phenomenon of extinction of $u$, alias blow-up of $1/u$, has been addressed by Fila and Stuke~\cite{FilaStuke14,FilaStuke16}. As a very simple paradigm for blow-up itself, consider the one-dimensional semilinear heat equation
\begin{equation}\label{C8:eq:nonlinear_heat}
u_t=u_{xx}+u^2
\end{equation}
on a bounded or unbounded interval.\ For real $u,x$ and \emph{real time} $t$, this equation has already been well-studied in the literature~\cite{Quittner07}. The solution is analytic in time and admits a local analytic continuation into the \emph{complex time plane}. For Dirichlet boundary conditions there exists a unique positive equilibrium $u^*$ which is one-dimensionally unstable. Solutions in the one-dimensional analytic unstable manifold of $u^*$ possess a real solution with complete blow-up in finite real time. Motivated by analyticity, Stuke studied extensions to complex, rather than real, time $t$. Of course this entails complex $u$. Following an idea of Masuda~\cite{Masuda84} and Guo at al.~\cite{Guo12}, Stuke attempts to bypass real time blow-up at $t=T$ by an excursion into the complex plane; see Fig.~\ref{fig:complex_time_domain}. Stuke proves that solutions on the complex unstable manifold of $u^*$ stay bounded on time paths $\Gamma_p^\pm$ parallel to the real axis. The solutions converge to zero: they are heteroclinic from $u^*$ to zero, along $\Gamma_p^\pm$, and locally foliate the nonsingular part of the complex one-dimensional unstable manifold $W_*^u$ of $u^*$. After some real time $t=T+\delta>T $, at the latest, the solutions $\Gamma^{\pm}_c$ extend back to the real axis from, both, the upper and the lower half-plane. Complete real blow-up, however, prevents these conjugate complex analytic branches to coincide, at real time overlap $t>T+\delta$. \\

Traveling pulse waves $u=U(x-ct)$ with wave speed $c$, and self-similar solutions, open a second approach to blow-up in complex time. Specifically, they correspond to ODE solutions which are homoclinic to zero. The standing wave $u(t,x)=U(x)=-6/x^2$ with wave speed $c=0$, for example, is an explicit solution of~\eqref{C8:eq:nonlinear_heat}. In general, homoclinic traveling waves $u=U(\psi)$, $\psi=x-ct $, possess local expansions in $\psi$ and $\log\psi$. \\

For the branched real blow-up on the complex unstable manifold $W^u_*$ of $u^*$, matters are not quite that simple. Indeed a "center manifold", rather than a fully hyperbolic structure, appears in rescaled variables. This heralds essential singularities beyond $\log t$, and sectors with exponentially flat correction terms. For further details see the dissertation~\cite{Stuke16}. \\

\subsection{The tumbling universe}

Arguably, the big-bang is the largest imaginable blow-up. For us it is the initial singularity of a purely gravitational cosmological model. We briefly describe the general setting, and the exact reduction to ODE models of Bianchi class. In the next section we present our own results. \\

Cosmological models of general relativity are solutions of the full Einstein equations
\begin{equation}\label{C8:eq:Einstein}
\mathrm{Ric}(g) - \tfrac{1}{2} \mathrm{Scal}(g) g \;\;=\;\; T\,,
\end{equation}
which relate the geometry of spacetime -- the curvature of a 4-dimensional Lorentzian metric $g$ -- to the matter content encoded in the stress-energy tensor $T$, with or without a cosmological constant. The Einstein PDEs are of hyperbolic type and are usually coupled to kinetic equations of matter. A complete understanding of the full PDE system~\eqref{C8:eq:Einstein} is currently out of reach. Therefore, additional symmetries are often assumed to discuss special solutions, and their perturbations. \\

The simplest model -- and core of the standard cosmological model still in use -- is the Friedmann model of spatially homogeneous and isotropic space-time. This assumption of a 6-dimensional symmetry group allows a reduction of~\eqref{C8:eq:Einstein} to a single scalar ODE that determines the expansion rate of the universe. This expansion rate can be compared with measurements of the Hubble constant and with the consequences of large-scale thermodynamics of the matter part. \\

The second simplest class are \emph{Bianchi models} of spatially homogeneous but anisotropic spacetime. In other words, the spacetime is assumed to be foliated into spatial hypersurfaces given by the orbits of a three-dimensional symmetry group. Belinskii, Khalatnikov, and Lifschitz (BKL) described the dynamics of cosmological models near the big-bang singularity as vacuum dominated, local, and oscillatory~\cite{BelinskiiKhalatnikovLifshitz1970-OscillatoryApproach}. This suggests to approximate the early universe by the dynamics of homogeneous, but not necessarily isotropic, vacuum spacetime. \\

The Einstein system~\eqref{C8:eq:Einstein} with a non-tilted, perfect-fluid matter model can be reduced to a 5-dimensional ODE system in expansion-reduced variables,
\begin{equation}\label{C8:eq:Bianchi}
\begin{array}{rcl}
 N_1' &=& ( q - 4 \Sigma_+ ) N_1\,,
\\ N_2' &=& ( q + 2 \Sigma_+ + 2\sqrt{3}\Sigma_- ) N_2\,,
\\ N_3' &=& ( q + 2 \Sigma_+ - 2\sqrt{3}\Sigma_- ) N_3\,,
\\ \Sigma_\pm' &=& (q-2) \Sigma_\pm - 3 S_\pm\,,
\end{array}
\end{equation}
with the abbreviations
\begin{equation}\label{C8:eq:BianchiAbbrv}
\begin{array}{rcl}
 S_+ &=& \tfrac{1}{2}\left(
 \left( N_2 - N_3 \right)^2
 - N_1 \left( 2 N_1 - N_2 - N_3 \right)
 \right)\,,
\\ S_- &=& \tfrac{1}{2}\sqrt{3}
 \left( N_3 - N_2 \right) \left( N_1 - N_2 - N_3 \right)\,,
\\ q &=& 2 \left( \Sigma_+^2 + \Sigma_-^2 \right)
 + \tfrac{1}{2}(3\gamma-2)\Omega\,,
\\ \Omega &=& 1 - \Sigma_+^2 - \Sigma_-^2 - K\,,
\\ K &=& \frac{3}{4} \Big( N_1^2 + N_2^2 + N_3^2
 - 2\left( N_1 N_2 + N_2 N_3 + N_3 N_1 \right)
 \Big)\,.
\end{array}
\end{equation}
This system is due to Wainwright and Hsu~\cite{WainwrightHsu1989-BianchiA}. The initial big-bang singularity is approached in the rescaled time limit $t\to-\infty$. The variables $N_k$ describe the curvature of the spatial hypersurfaces. Their signs determine the Lie-algebra type of the associated spatial symmetry imposed by the homogeneity assumption. Due to Bianchi's classification of three-dimensional Lie algebras -- the tangent spaces to the assumed symmetry group -- these models are called Bianchi models, although they have been introduced by G{\"o}del~\cite{Goedel1949-Cosmology} and Taub~\cite{Taub1951-SpaceTimes}. The shear variables $\Sigma_\pm$ relate to the rescaled eigenvalues of the second fundamental form of the spatial hypersurfaces. The matter density $\Omega$ is nonnegative, and the vacuum boundary $\Omega=0$ is invariant. The coefficient $\gamma<2$ determines the equation of state of the perfect fluid, e.g.\ $\gamma=4/3$ for radiation and $\gamma=1$ for dust. See also~\cite{WainwrightEllis1997-Cosmology} for further details on this dynamics approach to cosmology, and~\cite{HeinzleUggla2009-Mixmaster} for a survey on Bianchi models and open questions. \\

\begin{figure}
\centering
 \includegraphics[angle=0,scale=1.0]{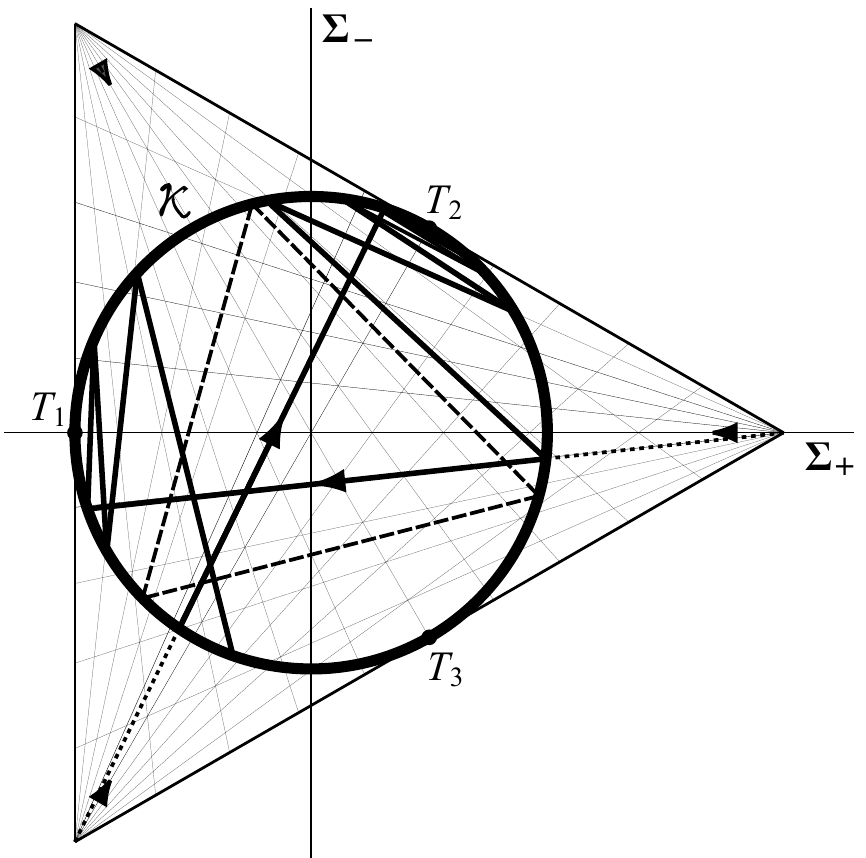}
\caption{\emph{The Kasner map and the three Taub points $T_1, T_2, T_3$ given by $N_1=N_2=N_3=0$, $(\Sigma_+ + i \Sigma_-)^3=-1$.}}
\label{C8:fig:Kasner} 
\end{figure}

The most prominent features of system~\eqref{C8:eq:Bianchi} are the \emph{Kasner circle} $\mathcal{K}$ of equilibria,
\begin{equation}\label{C8:eq:Kasner}
\mathcal{K} \;=\; \{\; \Sigma_+^2+\Sigma_-^2=1, \; N_1=N_2=N_3=0 \;\}\,,
\end{equation}
and the invariant \emph{Kasner caps}
\begin{equation}\label{C8:eq:Caps}
\mathcal{H}_k^\pm \;=\; \{\; \Sigma_+^2+\Sigma_-^2=1-N_k^2, \;
 \pm N_k > 0, \; N_{k+1}=N_{k-1}=0, \; k \;\mathrm{mod}\; 3 \;\}\,.
\end{equation}
\begin{figure}
\centering
\includegraphics[angle=0,scale=1.1]{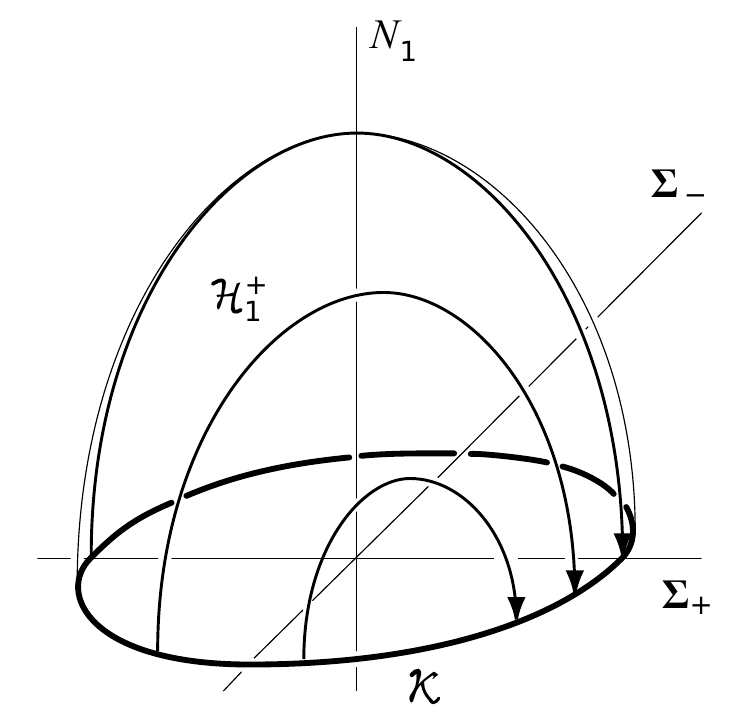}
\caption{\emph{Heteroclinic cap of vacuum Bianchi II solutions to the Kasner circle $\mathcal{K}$. Remaining caps are given by $2\pi/3$ rotations over $2\pi/3$ in $\Sigma_\pm$ and cyclic permutations of $\{N_1, N_2, N_3\}$.}}
\label{C8:fig:Caps} 
\end{figure}
The Kasner caps are filled with heteroclinic orbits between equilibria on $\mathcal{K}$; see Fig.~\ref{C8:fig:Caps}. Equilibria of~\eqref{C8:eq:Bianchi} indicate self-similar blow-up in the big-bang limit which occurs at rescaled time $t\to -\infty$. The projections of the heteroclinic orbits to the $\Sigma$-plane lie on straight lines through the corners of a circumscribed equilateral triangle; see Fig.~\ref{C8:fig:Kasner}. In reversed physical time, these sequences define the \emph{Kasner map} as an expansive map from the Kasner circle onto its double cover; see Fig.~\ref{C8:fig:Kasner} and Section~\ref{subsec:kasnerdynamics}. \\

Concatenation of the heteroclinic orbits in the Kasner caps, i.e.\ iteration of the Kasner map, defines a subshift of finite type. Misner~\cite{Misner1969-Mixmaster} introduced the name \emph{Mixmaster} for the heuristic picture of generic trajectories following these formal heteroclinic sequences close to the big-bang singularity: the universe tumbles between various close-to self-similar solutions given by the various Kasner equilibria. \\

A crucial question is the \emph{locality} of the BKL/Mixmaster picture. In the (non-vacuum) Friedmann model, the particle horizon (backward light cone) remains finite in expansion reduced variables, close to the big-bang singularity. This bounded particle horizon also bounds the domain of dependence by particle interaction, and therefore represents locality. In the Bianchi model~\eqref{C8:eq:Bianchi}, the particle horizon corresponds to the \emph{Mixmaster integral}
\begin{equation}\label{C8:eq:MixIntegral}
I \;\;=\;\; \int_{-\infty}^0(\sqrt{N_1N_2}+\sqrt{N_2N_3}+\sqrt{N_1N_3}) \mathrm{d}t\,,
\end{equation}
see~\cite{HeinzleRingstroem2009-FutureBianchiVI}. Locality holds true, if and only if this integral remains finite. Obvious counterexamples are the
solutions on the line of Taub equilibria
\begin{equation}
N_1=0\ne N_2=N_3\,,\qquad \Sigma_-=0\,,\quad\Sigma_+=-1\,,\qquad \Omega=0\,,
\end{equation}
where the Mixmaster integral is $+\infty$. It is currently not known, whether nontrivial counterexamples with $I=\infty$ exist. The first nontrivial positive examples with $I<\infty$ are quite recent~\cite{LiebscherHaerterichWebsterGeorgi2011-BianchiA}. More precisely, the longstanding BKL/Mixmaster conjecture therefore claims that $I<\infty$ for generic initial conditions $N_j(0)$, $\Sigma_{\pm}(0)$ in the respective Bianchi class~\cite{BelinskiiKhalatnikovLifshitz1970-OscillatoryApproach,HeinzleRingstroem2009-FutureBianchiVI}. Here \emph{genericity} can be either understood as \emph{residual}, in the topological sense of Baire category, or as \emph{almost everywhere} with respect to Lebesgue measure. In the Bianchi type IX of class A, where all $N_j>0$, Ringstr\"om~\cite{Ringstroem2001-BianchiIX} has shown Baire and Lebesgue generic convergence to the Kasner caps. Of course, the Mixmaster integral $I$ vanishes on the Kasner caps. The mere generic convergence to the Kasner caps, however, does not imply generic finiteness of $I$ if this convergence is too slow.

\subsection{Kasner dynamics}\label{subsec:kasnerdynamics}

In the vacuum model~\eqref{C8:eq:Bianchi}, $\Omega=0$, each non-Taub Kasner equilibrium possesses one trivial eigenvalue zero of the linearization in the direction of the Kasner circle and three nontrivial eigenvalues in the remaining directions of the three heteroclinic caps. At the Taub points, two nontrivial eigenvalues become zero and a \emph{bifurcation without parameters} occurs, in the sense of the habilitation thesis~\cite{Liebscher15-Bifurcation} of Stefan~Liebscher. \\

To clarify the relation of generic trajectories to formal heteroclinic sequences, we ask for the set of points converging to the heteroclinic sequence attached to a given Kasner equilibrium. This is a difficult question because each heteroclinic step, by itself, already consumes infinite rescaled time $t$. We might hope for an invariant codimension-one foliation transverse to the Kasner circle, which in particular contains the local strong unstable and strong stable invariant manifolds of each Kasner equilibrium. Each fiber would then become the stable set of the Kasner equilibrium, in backwards time $t \to -\infty$, as the big-bang singularity is approached. \\

However, there are topological obstructions: the sets of initial conditions converging to the Taub points, backwards, were classified in~\cite{Ringstroem2001-BianchiIX} and form manifolds of codimension two. Thus, if the heteroclinic sequence of a Kasner equilibrium terminates at a Taub point, this equilibrium cannot possess an attached stable fiber of codimension one. In particular, any transverse backwards stable foliation has to miss at least a dense subset of the Kasner circle. \\

Stable foliations along Kasner sequences away from Taub points have first been addressed by Liebscher et al.\ in~\cite{LiebscherHaerterichWebsterGeorgi2011-BianchiA}, as follows. Consider a Kasner equilibrium such that its attached (unique, backwards) heteroclinic sequence does not accumulate at a Taub point. Then this equilibrium possesses a unique, locally invariant, codimension-one, stable Lipschitz leaf. The leaves of equilibria in the same heteroclinic sequence form a backwards invariant continuous foliation. The distance in which the trajectories pass the equilibria in this foliation, tends to zero exponentially. \\

The proof basically consists of three parts. First, a suitable eigenvalue property guarantees an arbitrarily strong contraction, transversely to the Kasner circle, and an arbitrarily small drift along the Kasner circle, locally. Second, the global passage along a heteroclinic orbit in the Kasner caps is modelled by a diffeomorphism close to the Kasner map. This provides a bounded expansion along the Kasner circle, away from Taub points. Third, the combination of both maps yields a hyperbolic map between suitably chosen cross sections near the Kasner circle: contraction transverse to the Kasner circle is inherited from the local passage, expansion along the Kasner circle is inherited from the global passage. The stable leaf is then obtained as the fixed point of a graph transform, for the above construction. \\

In~\cite{LiebscherRendalTchapnda2012-BianchiMaxwell}, Liebscher et al.\ generalize this approach to weaker eigenvalue conditions. This applies to non-vacuum Bianchi systems with arbitrary ideal fluids, and to Einstein--Maxwell equations of Bianchi type VI${}_0$. In the Einstein--Maxwell system one of the heteroclinic caps is no longer induced by the geometry but rather by the dynamics of the magnetic field. \\

An alternative approach in~\cite{Beguin2010-BianchiAsymptotics}, and the dissertation~\cite{Buchner2013-Dissertation} by Buchner, utilize Takens linearization~\cite{Takes1971-Linearization} to study the local passage. Indeed, after linearization, contraction and absence of drift follow trivially. However, the linearization requires strong spectral non-resonance conditions. Non-resonance imposes severe limitations on this approach and on possible generalizations. In~\cite{Beguin2010-BianchiAsymptotics}, for example, all formal sequences which are periodic or accumulate to periodic sequences (or the Taub points) are excluded to avoid resonances. In~\cite{Buchner2013-Dissertation}, Buchner has studied the non-resonance conditions more carefully, and certain periodic sequences are now admissible. His results on local passages extend to Bianchi class-B systems. \\

The techniques described above require to avoid the Taub points. This is their main restriction. In fact, the Kasner equilibria which, under iteration of the Kasner map, do not accumulate at Taub points form an (albeit uncountable) meager set of measure zero. They are not generic, in any sense, but possess a finite Mixmaster integral~\eqref{C8:eq:MixIntegral} and thus a finite particle horizon. \\

Brehm, in his dissertation~\cite{Brehm16}, is the first to answer the BKL/Mixmaster question in the Lebesgue sense. This problem had remained open for more than four decades. Brehm obtains two new results in the study of the Wainwright--Hsu system~\eqref{C8:eq:Bianchi},~\eqref{C8:eq:BianchiAbbrv} for the vacuum case $\Omega=0$.The first result extends Ringstr\"om's attractor result to Bianchi type \textsc{VIII} vacuum solutions. The second result asserts finiteness $I<\infty$ of the Mixmaster integral~\eqref{C8:eq:MixIntegral} for Lebesgue almost all initial conditions. \\

Consider the Bianchi type \textsc{IX} first, where all $N_j>0$ initially (and then for all times). The \emph{Taub space} $\{N_2=N_3,\, \Sigma_-=0\}$ is time invariant, together with its two rotated ``cousins'' given by cyclic permutations of the indices $\{1,2,3\}$. Ringstr\"om~\cite{Ringstroem2001-BianchiIX} has shown convergence to the Kasner caps, for initial conditions outside the Taub spaces, in the time-rescaled big-bang limit $t\to-\infty$. The proof proceeded by averaging over rotations around the invariant Taub-spaces, which are known explicitly. \\

The Bianchi type \textsc{VIII} is defined by $N_2,N_3>0>N_1$ (up to index permutations). The above Taub-space remains invariant, but its two cousins with $N_1=N_2$ and $N_1=N_3$ are now incompatible with $N_2,N_3>0>N_1$ and have disappeared. Even worse: all solutions with $N_3=0$ and $N_2>0>N_1$ converge to the Kasner-cicle $\mathcal K$ in \emph{both} time-directions $t\to\pm\infty$ and fail to lie on the Kasner-caps. This poses a serious obstacle to the attractivity of the caps. For these reasons, attractor results have remained elusive in the Bianchi \textsc{VIII} setting. \\

For both Bianchi \textsc{VIII} ($N_2,N_3>0>N_1$) and Bianchi \textsc{IX} ($N_1,N_2,N_3>0$), the attractor problem for Bianchi class A vacuum ($\Omega=0$) has been settled by Brehm in his dissertation~\cite{Brehm16} as follows. The dynamics in the big-bang limit $t\to-\infty$ falls into one of the following three mutually exclusive classes:
\begin{itemize}
\item [$(a)$] The solution converges to the three Kasner caps and has at least one $\alpha$-limit point on the Kasner-circle, which is distinct from the three Taub-points $T_j$ of Fig.~\ref{C8:fig:Kasner} and their antipodal points $Q_j=-T_j$.
\item [$(b)$] The solution is contained in the invariant Taub subspace $\{N_2=N_3,\, \Sigma_-=0\}$, or, in the case of Bianchi \textsc{IX}, one of its two cousins.
\item[$(c_2)$] The solution has exactly two $\alpha$-limit points on the Kasner-circle, which are $T_2$ and $Q_2=-T_2$. It has $\limsup |N_1N_3|>0 = \liminf|N_1N_3|$, for $t\to-\infty$, and $\lim|N_1N_2|=\lim|N_2N_3|=0$.
\item[$(c_3)$] The analogon of $(c_2)$, with the indices $2$ and $3$ interchanged.
\end{itemize}
Case $(a)$ applies for an open set of initial conditions which are a neighborhood of the three caps minus the three Taub spaces. Furthermore, for any solution in case $(a)$, the following non-Mixmaster integral $J$ remains finite:
\begin{equation}\label{C8:eq:brehm:weak-integral}
J=\int_{-\infty}^0 \left(|N_1N_2|+|N_2N_3|+|N_3N_1|\right)\mathrm{d} t < \infty\,.
\end{equation}
The cases $(c_2)$ and $(c_3)$ are impossible in the case of Bianchi \textsc{IX}; whether they can occur at all is currently unknown. Taken together with suitable measure-theoretic estimates, one can see that case $(a)$ applies in a set of initial conditions, both of fat Baire category and of full Lebesgue measure. In other words $(b),(c_2),(c_3)$ occur for a meager set of Lebesgue-measure zero, only. The estimate~\eqref{C8:eq:brehm:weak-integral} for the non-Mixmaster integral $J$ is novel also in the Bianchi \textsc{IX} case. It is insufficient to bound the Mixmaster integral $I$ of~\eqref{C8:eq:MixIntegral}, but nevertheless promising: from $\lim_{t\to-\infty}|N_iN_j|= 0$, to $\int_{-\infty}^0|N_iN_j|\mathrm{d} t <\infty$, and finally towards the goal $\int^0_{-\infty}\sqrt{|N_iN_j|}\mathrm{d} t <\infty$. 
The proof in Bianchi \textsc{IX}, where only the $J<\infty$ estimate~\eqref{C8:eq:brehm:weak-integral} is novel compared to~\cite{Ringstroem2001-BianchiIX}, proceeds along the lines of Ringstr\"om, via a refined quantitative averaging over rotations around the three invariant Taub subspaces combined with a more careful gluing of the different averaging regimes. 
In the Bianchi \textsc{VIII} case, the two missing invariant Taub subspaces are replaced by $\{|N_1|=|N_3|,\, \Sigma_-=\sqrt{3}\Sigma_+\}$ and its index-exchanged cousin. These two replacement Taub subspaces are invariant only up to leading order; this suffices for the averaging arguments and yields the trichotomy $(a)$, $(b)$, $(c)$. 
The second main result in~\cite{Brehm16} establishes boundedness $I<\infty$ of the Mixmaster integral~\eqref{C8:eq:MixIntegral}, for Lebesgue almost every initial condition in Bianchi \textsc{VIII} and \textsc{IX} vacuum. In other words, particle horizons form towards the singularity and the blow-up is local in the sense of BKL. \\

\begin{figure}
\begin{center}
\includegraphics[scale=0.8]{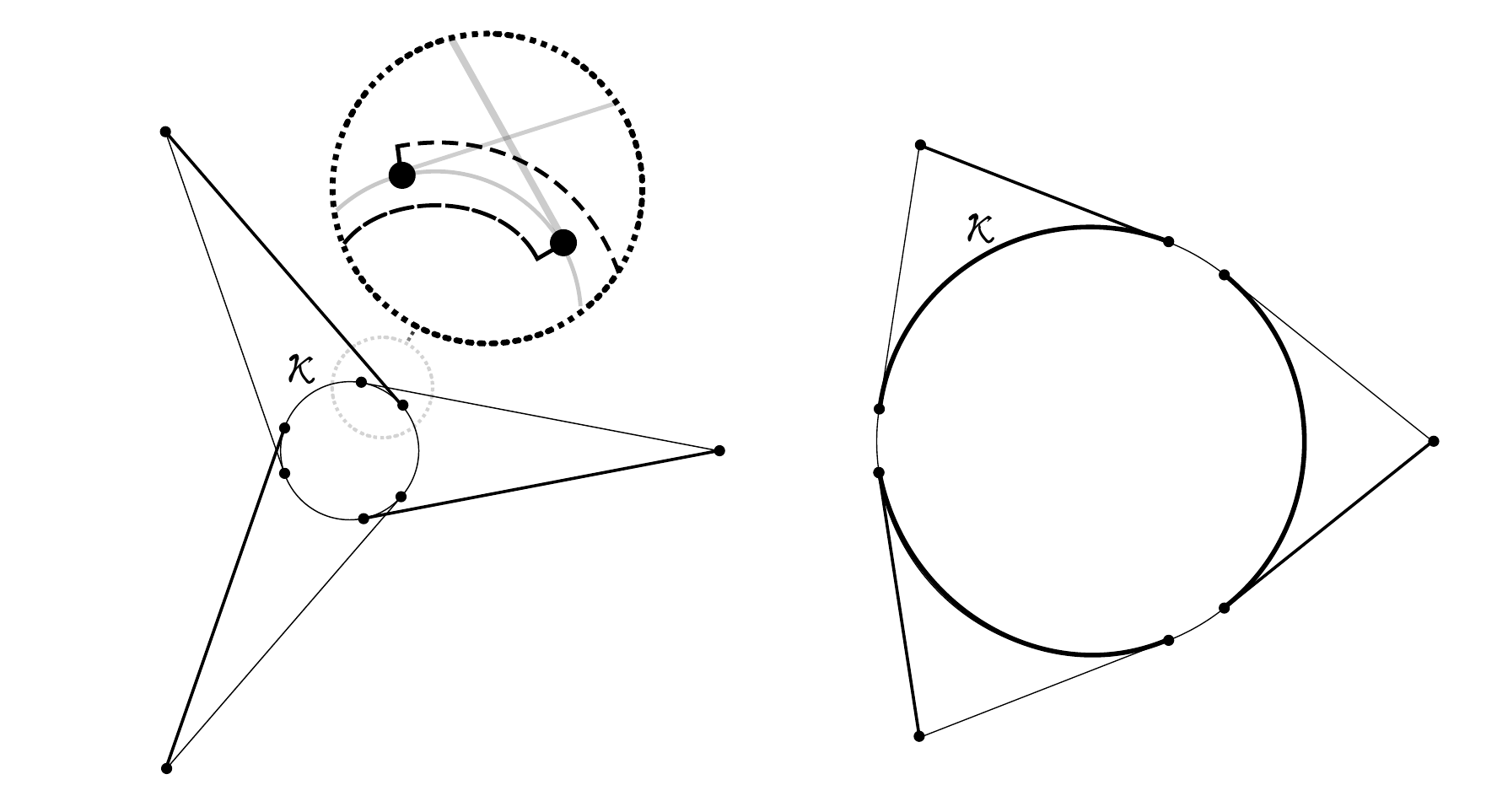}
\end{center}
\caption{\emph{The Kasner circle of equilibria, the three expanding arcs and the three points of emanation for the HL models. Left $0<v<1/2$; the dashed circle zooms into one of the three overlapping regions, where the overlapping sectors are indicated by dashed lines. Right $1/2<v<1$.}}
\label{HLKasnerdyn}
\end{figure}

The above considerations reside within the general relativity (GR) framework of a purely gravitational Einstein universe. A larger framework is provided by the Ho\v{r}ava-Lifshitz models (HL). HL gravity has been proposed as a candidate for the theory of quantum gravity; see~\cite{Muko} and the references there. The motivation of Hell to study those models in~\cite{HellUggla}, with Uggla, is the following. HL gravity breaks relativistic first principles and introduces anisotropic scalings between space and time. How do such modifications of the first principles affect the initial singularity and the structure of the BKL conjecture? As an attempt to answer this question, Hell and Uggla address the question of the dynamics on its first level of complexity, i.e.\ the discrete dynamics contained in the heteroclinic chains arising in HL. Similarly to the Bianchi models for general relativity, it is possible to formulate HL models as a system of ODEs with constraints that admit a circle of equilibria, which we call the Kasner circle by analogy. It corresponds to the Bianchi type I model for HL. On the next level of the symmetry hierarchy, Bianchi type II, we find three invariant hemispheres of heteroclinic orbits that intersect only at their boundary, the Kasner circle, for HL and GR alike. On each hemisphere, the projections of the heteroclinic orbits to the plane of the Kasner circle emanate from a single point, as a straightforward calculation shows\cite{HellUggla}. Those three points of emanation form an equilateral triangle whose center coincide with the center of the circle of equilibria: for the HL parameter $v=1/2$ corresponding to GR, the Kasner circle of equilibria is inscribed in the triangle, while this is not the case any more when the parameter changes to $v\neq 1/2$; see Fig.~\ref{HLKasnerdyn}. Concatenation of heteroclinic orbits on different hemispheres form heteroclinic chains, just as they did in the GR case. However, the discrete dynamics on the Kasner circle induced by those heteroclinic chains change. \\

The BKL conjecture for GR describes the asymptotic behavior towards the big-bang singularity. This conjecture can be seen as an extrapolation from the dynamics of the heteroclinic chains to the continuous dynamics near the invariant hemispheres. Studying the discrete dynamics for the HL Bianchi models reveals how the BKL conjecture should be modified in this context. \\

Let us focus on the discrete dynamics induced by the heteroclinic chains, which map each equilibrium of the Kasner circle to the end of the heteroclinic leaving from it. Each point of emanation expands an arc of the Kasner circle enclosed by the two tangents to the circle through this point, by this HL Kasner map. Varying $v$ amounts to changing the distance between the points of emanation of the heteroclinics and the circle of equilibria. This changes the size of the three expanding arcs and the factor of expansion; see Fig.~\ref{HLKasnerdyn}. Two very different behaviors appear, depending on whether $v$ is below or above $1/2$. \\

Let us first discuss the case $1/2<v<1$. As can be seen in Fig.~\ref{HLKasnerdyn}, right, the three expanding arcs are separated by three arcs of stable fix points at which heteroclinic chains end, due to the lack of heteroclinic orbits leaving from there. As proved in~\cite{HellUggla}, generic initial conditions (both in the sense of Baire and Lebesgue) on the Kasner circle will lead to finite heteroclinic chains that end in a stable arc. Hence an analogue to the BKL conjecture should state that the asymptotic behavior of the system is, in general, convergence to a single equilibrium in a stable arc. This behavior is very different from the GR case where oscillatory chaos is present in the asymptotic behavior towards the big-bang singularity. However, by~\cite{HellUggla}, not all initial conditions on the Kasner circle lead to finite heteroclinic chains. An invariant Cantor set of measure zero leads to heteroclinic chains that never enter the stable arcs, and the discrete dynamics on this Cantor set is chaotic. The proofs of these facts use symbolic dynamics and equivalence to a subshift of finite type, known to be chaotic. Hence chaos has not completely vanished from the system for $1/2<v<1$, but it is not generic. \\

Let us now discuss the opposite case $0<v<1/2$. We see from Fig.~\ref{HLKasnerdyn}, left, how the three expanding arcs now overlap, pairwise, while no stable arc is present. Therefore the continuation of a heteroclinic chain that hits the overlap region is not unique: there are two heteroclinics leaving from there. In other words, the map describing the discrete dynamics is not well-defined and the concept of iterated function systems (IFS) comes into play -- details can be found in~\cite{HellHL} by Hell. The IFS we consider is a collection of maps whose domains of definition cover the Kasner circle. The notion of chaotic discrete dynamical systems is generalized to IFS using the Hausdorff distance between sets. Iterations of a well-defined map are said to be chaotic if three conditions are fulfilled: sensitivity to initial conditions, topological mixing, and density of periodic orbits. In the case of an IFS, a single initial condition generates several trajectories. Therefore each iteration of an IFS is a set instead of a single state. The three above conditions for chaos translate to IFS when distance between iterations is translated to Hausdorff distance between sets of possible iterations. The IFS arising for the HL Bianchi models with $0<v<1/2$ are chaotic in this sense, mostly due to the expansion property. Roughly speaking, this means that there are chaotic realizations of the IFS. This property only concerns the heteroclinic chains. The nearby flow will be investigated elsewhere. Even if the HL models with $0<v<1/2$ are chaotic, like the GR case $v=1/2$, their chaotic structure is quite different. The concept of a BKL era (i.e.\ phases of oscillations between two neighboring expanding arcs) and the growth of their length along heteroclinic chains shape the chaos in the GR case, via the continued fraction expansion. In the HL for $0<v<1/2$, BKL erae do not play such a role: their lengths are bounded above, and the fine structure associated to their rate of growth disappears. The value $v=1/2$ associated to GR can be seen as a critical transitional value where the dynamics bifurcates from nongeneric chaos to generic chaos. The importance of erae is a feature that does not survive below $v=1/2$. Hence the perturbation of GR along the one-parameter family of HL will lead to a breakdown of the BKL conjecture. \\

This concludes our extended survey of some global dynamics aspects, from geometric and algebraic aspects of parabolic global attractors and their blow-up to black hole data and exact chaotic solutions of the Einstein equations in the setting of Bianchi cosmology. Not only many of the particular problems but, much broader, our overall perspective on these problems we owe, deeply, to the unique framework provided by the Sonderforschungsbereich 647 ``Space -- Time -- Matter'' of the Deutsche Forschungsgemeinschaft.

\bibliographystyle{amsplain}
\bibliography{Fiedler_etal_SFB647}
\end{document}